\documentclass[12pt,titlepage]{article}
\usepackage{graphicx,color}
\usepackage{amsmath,amssymb,amsthm}
\input epsf

\newtheorem{thm}{Theorem}
\newtheorem{lemma}[thm]{Lemma}
\newtheorem{cor}[thm]{Corollary}
\newtheorem{problem}{Problem}

\newtheorem{prop}[thm]{Proposition}

\newtheorem{conj}{Conjecture}
\newtheorem{defi}{Definition}

\def\beq{\begin{equation}}\def\eeq{\end{equation}}
\def\beqn{\begin{eqnarray}}\def\eeqn{\end{eqnarray}}

\def\qed{\ifhmode\unskip\nobreak\fi\quad\ifmmode\Box\else$\Box$\fi}

\newcommand\SG{{\rm SG}}
\newcommand\reals{{\mathbb{R}}}
\title{On $4$-chromatic Schrijver graphs: their structure,
  non-$3$-colorability, and critical edges}

\author{{\bf G\'abor Simonyi}$^{a,b,}$\thanks{Research partially supported by the
    National Research, Development and Innovation Office (NKFIH)
    grants K--116769, K--120706 and by the National Research, Development and Innovation Fund (TUDFO/51757/2019-ITM, Thematic Excellence Program).}
    $\qquad$ {\bf G\'abor Tardos}$^{a,c,}$\thanks{Research partially supported by the Cryptography ``Lend\"ulet'' project of the Hungarian Academy of Sciences and by the National Research, Development and Innovation Office (NKFIH) grants K--116769 and SSN-117879.}\\ \\
$^a$Alfr\'ed R\'enyi Institute\\ \\
$^b$Department of Computer Science and Information Theory,\\ Budapest University of Technology and Economics\\ \\
$^c$Department of Mathematics,\\
Central European University\\ \\
 {\tt simonyi@renyi.hu} \ \ \ {\tt tardos@renyi.hu}}
\date{}

\begin{document}
\maketitle
\begin{abstract}

We give an elementary proof for the non-$3$-colorability of $4$-chromatic
Schrijver graphs thus providing such a proof also for $4$-chromatic Kneser
graphs. To this end we use a complete description of the structure of
$4$-chromatic Schrijver graphs that was already given by Braun in
\cite{BBraun1, BBraun2} and even earlier in an unpublished manuscript by Li
\cite{Li}. We also address
connections to surface
quadrangulations. In particular, we show that a spanning subgraph of
$4$-chromatic Schrijver graphs quadrangulates the Klein bottle, while another
spanning subgraph quadrangulates the projective plane. The latter is a
special case of a result by Kaiser and Stehl\'{\i}k \cite{KS2}. We
characterize the color-critical edges of $4$-chromatic Schrijver graphs and also present preliminary results toward the characterization of color-critical edges in Schrijver graphs of higher chromatic number. Finally, we show that (apart from two cases of small parameters) the subgraphs we present that quadrangulate the Klein bottle are edge-color-critical.
The analogous result for the subgraphs quadrangulating the projective plane is an immediate consequence of earlier results by Gimbel and Thomassen and was already noted by Kaiser and Stehl\'{\i}k.

\bigskip
\bigskip
\par\noindent

\end{abstract}

\section{Introduction}
\message{Introduction}

Lov\'asz's celebrated proof of Kneser's conjecture \cite{LLKn} is based on the Borsuk-Ulam theorem from algebraic topology and its appearance triggered a lot of activities in applying algebraic topological tools, in particular, the Borsuk-Ulam theorem in combinatorics. For several examples of this the reader is referred to Ji\v{r}\'{\i} Matou\v{s}ek's excellent book \cite{Matbook}. Naturally, it was always a question, whether the topological tools are really necessary to prove Kneser's conjecture, or a purely combinatorial argument could substitute it. Matou\v{s}ek himself gave such a proof in \cite{Matcombpr} that Ziegler \cite{ZieglInv} generalized later to make it apply to generalized versions of the Lov\'asz-Kneser theorem. Among others, Ziegler also gave a combinatorial proof of Schrijver's theorem \cite{Sch78}. These proofs, as Ziegler puts it, ``are combinatorial (`elementary') in the sense that they do not rely on topological concepts (such as continuous maps, simplicial approximation, homology) or results.'' Admittedly, however, all these proofs are inspired by topology and (quoting Ziegler again) ``it is neither desirable nor practical to eliminate this background intuition.''

One can still ask, however, whether the topological background is also
necessary for such special cases of these theorems, where the parameters
involved are small. Kneser graphs depend on two positive integer parameters,
$n$ and $k$ (satisfying $n\ge 2k$; see Definition~\ref{defi:kg} below) and the
Lov\'asz-Kneser theorem states that their chromatic number is $n-2k+2$. The
difficult part of the proof is the lower bound part, that is proving that this
many colors are needed for a proper coloring. This is completely trivial,
however, when this number is $2$, and is also quite easy when it is $3$, as it
only needs the demonstration of an odd cycle in the graph. (In case of Schrijver graphs, their $3$-chromatic
versions are just the odd cycles themselves.) The $4$-chromatic case, however, is already non-trivial, and this is what we focus on in this paper.

To be able to say more let us give the necessary definitions. We start with standard notation.
For $n\ge1$, we write $[n]$ to denote the set of positive integers up to
$n$. For $k\ge0$ and a set $S$ we write ${S\choose k}$ for the set of size $k$
subsets of $S$. In this paper we deal with simple graphs only. All the graphs mentioned are finite except for the Borsuk graphs used in Section~\ref{otegy}. If $G$ is a
graph we write $V(G)$ for its set of vertices and $E(G)$ for the set of
edges. An edge connecting the vertices $a$ and $b$ is denoted by $ab$ (or
equivalently by $ba$). A map $f$ from the vertex set of a graph is called a \emph{coloring} of the graph and the value $f(v)$ is called the \emph{color} of the vertex $v$. We call an edge \emph{monochromatic} if it connects two vertices of the same color and we call a coloring \emph{proper} if no edge is monochromatic. We write $\chi(G)$ for the \emph{chromatic number} of
the graph $G$, that is, the smallest number of colors in a proper coloring of $G$. In order to be able to name the vertices of paths and cycles we
use $P_m$ and $C_m$ for path or cycle specifically on the vertex set
$[m]$. Namely for $m\ge2$ integer we write $P_m$ for the path with vertex set
$[m]$ and with edges connecting consecutive numbers. For $m\ge3$ we write
$C_m$ for the cycle obtained from $P_m$ by adding the edge $m1$. For a graph $G$ and a subset $S$ of $V(G)$ we write $G[S]$ for the subgraph induced by $S$ and $G\setminus S=G[V(G)\setminus S]$. Similarly, $G\setminus S$ for $S\subseteq E(G)$ means the subgraph obtained from $G$ by deleting the edges in $S$. For graphs $G$ and $H$ a map
$f:V(G)\to V(H)$ is called a \emph{homomorphism} if for every edge $ab\in E(G)$ $f(a)f(b)$ is an edge in $H$. Clearly, if a homomorphism from $G$ to $H$ exists, then we have $\chi(G)\le\chi(H)$. We call a map $f:V(G)\to V(H)$ a \emph{cover} if for any vertex $v$ of $G$ $f$ maps the neighbors of $v$ in $G$ bijectively to the neighbors of $f(v)$ in $H$. We call such a cover a \emph{double cover} if every vertex in $H$ is the image of exactly two vertices of $G$.

\begin{defi}\label{defi:kg}
Let $n\ge 2k$ be positive integers.
The \emph{Kneser graph} ${\rm KG}(n,k)$ is defined as follows.
\smallskip

$V({\rm KG}(n,k))={[n]\choose k}$,

$E({\rm KG}(n,k))=\{ab: a,b\in {[n]\choose k}, a\cap b=\emptyset\}$.

The \emph{Schrijver graph} ${\rm SG}(n,k)$ is defined as the subgraph of ${\rm KG}(n,k)$ induced by the vertices of the Kneser graph that are independent sets in the cycle $C_n$.
\end{defi}
\medskip

Kneser \cite{Kne} in 1955 observed (using different terminology) that the chromatic number $\chi({\rm KG}(n,k))$ is at most $n-2k+2$ and conjectured that this inequality is sharp. This was proved by Lov\'asz \cite{LLKn} in 1978.

\begin{thm} (Lov\'asz-Kneser theorem)
$$\chi({\rm KG}(n,k))=n-2k+2.$$
\end{thm}

Soon afterwards, Schrijver generalized (actually B\'ar\'any's simplified version \cite{Barany} of) the proof and showed that the induced subgraphs now named after him have the same chromatic number. He also observed that these graphs are already vertex-color-critical.

\begin{thm} (Schrijver \cite{Sch78})
$$\chi({\rm SG}(n,k))=n-2k+2.$$

Furthermore, for every $a\in V({\rm SG}(n,k))$ we have $\chi({\rm SG}(n,k)\setminus\{a\})=n-2k+1,$.
\end{thm}

\medskip
In this paper we give an elementary proof of the statement that $\chi({\rm SG}(2k+2,k))\ge 4$. This proof is combinatorial in the sense that it can be followed without any knowledge in topology. Nevertheless, an expert can feel the topological connection even in this argument:
we will use a parameter called winding number that, though completely
elementary and can be understood without any topological background,
undeniably bears some flavour of topology. Our proof is very closely
related to Mohar and Seymour's proof in \cite{MS} of non-$3$-colorability of
certain surface quadrangulations, where the winding number also plays a
crucial role. In fact, $4$-chromatic Schrijver graphs have a close connection
to such quadrangulations that we will discuss in Section~\ref{sect:quadr}.

Our proof builds on an understanding of how $4$-chromatic Schrijver graphs
look like.
Such a structural description of ${\rm SG}(2k+2,k)$ is already implicit in Braun's proof of
his result in \cite{BBraun1} about the automorphism group of Schrijver graphs
and is made more  explicit in \cite{BBraun2}. It also appears in the even
earlier manuscript by Li \cite{Li}.
Nevertheless, in Section~\ref{sect:struct} we also present such a description
(in somewhat different and more general terms) for the sake of completeness.
Section~\ref{sect:proof} contains the actual proof of
non-$3$-colorability. Using the structural description from
Section~\ref{sect:struct} again, we show in Section~\ref{sect:quadr} that a spanning
subgraph of every $4$-chromatic Schrijver graph quadrangulates the Klein
bottle. This observation itself also provides a proof of
non-$3$-colorability if we use a theorem independently proved by Archdeacon,
Hutchinson, Nakamoto, Negami, Ota \cite{AHNNO} and Mohar and Seymour
\cite{MS}. With slight modification of the previous (and somewhat more natural) quadrangulation we will also see how another spanning subgraph of ${\rm SG}(2k+2,k)$
quadrangulates the projective plane. This gives a special case of the result
of Kaiser and Stehl\'{\i}k in \cite{KS2} and also implies non-$3$-colorability
by a result of Youngs \cite{Y}, which is generalized in \cite{AHNNO} and
\cite{MS}, and also in \cite{KS1} in a different direction.

As we have seen Schrijver proved that Schrijver graphs are
vertex-color-critical. They are not edge-color-critical, however. In
Section~\ref{sect:crit} we
characterize the color-critical edges of $4$-chromatic Schrijver graphs and also present preliminary results toward the characterization of color-critical edges in Schrijver graphs of higher chromatic number. In Section~\ref{critklein} we show that (except for very small values of the parameters) the special subgraphs of $\SG(2k+2,k)$ quadrangulating the Klein bottle are edge-color-critical. (The analogous result for the subgraphs quadrangulating
the projective plane is a straightforward consequence of a theorem in
\cite{GT} as noted in \cite{KS2}.)

As we wrote this paper we found that some ingredients of our argument (in
particular, the structure of graphs ${\rm SG}(2k+2,k)$ and the usefulness of
the winding number in similar  investigations) are already known, the relevant
references are given in the foregoing. Nevertheless, these ingredients seem
not to have ever appeared together the way we present them and therefore they
are given here in full detail.
Apart from summarizing these facts, the main novelties in this paper are the
way we present the structure of $4$-chromatic Schrijver graphs via the
more general family of graphs we call reduced drums (see Definition~\ref{defi:drum} in
Section~\ref{sect:struct}), their direct use in the proof of non-$3$-colorability
(Section~\ref{sect:proof}), our results on critical edges of Schrijver graphs
(Section~\ref{sect:crit}), and the edge-color-criticality of certain subgraphs
of four-chromatic Schrijver graphs quadrangulating the Klein bottle (Section~\ref{sect:quadr}). In the final short section we formulate a conjecture that would characterize critical edges in Schrijver graphs that are sufficiently large for their chromatic number.

\section{The structure of $4$-chromatic Schrijver graphs}\label{sect:struct}

To describe the structure of the graphs ${\rm SG}(2k+2,k)$, we will use the Cartesian product of graphs and also identification of vertices.

\begin{defi}\label{defi:Cartprod}
Let $F$ and $G$ be two graphs. Their \emph{Cartesian product} $F\Box G$ has vertex set $V(F\Box G)=V(F)\times V(G)$, while the edge set is defined by
$$E(F\Box G)=\{(a,c)(b,d): (a=b\in V(F)\ {\rm and}\ cd\in E(G))\ {\rm or}\ (ab\in E(F)\ {\rm and}\ c=d\in V(G))\}.$$
\end{defi}

The subgraphs of $F\Box G$ induced by the vertex sets $\{a\}\times V(G)$ for $a\in V(F)$ are called {\em layers}. Note that all layers are isomorphic to $G$. Symmetrically, the subgraphs induced by $V(F)\times \{c\}$ for $c\in V(G)$ are isomorphic to $F$, but we will not consider this type of subgraphs.

\begin{defi} \label{defi:ident}
Let $F$ be a graph and $a,b\in V(F)$. By \emph{identifying} vertices $a$ and $b$ in $F$ we obtain the graph with vertex set $(V(F)\setminus\{a,b\})\cup\{v_{a,b}\}$ in which two vertices are connected if neither is $v_{a,b}$ and they are connected in $F$ or one of them is $v_{a,b}$ and the other is a neighbor of either $a$ or $b$ in $F$. Clearly, we can also identify multiple pairs of vertices in a graph.
\end{defi}

The following structure will give a kind of skeleton of the graphs ${\rm SG}(2k+2,k)$, therefore we give them a name.

\begin{defi}\label{defi:drum}
The {\em drum} of height $h$ and perimeter $n$, denoted $D_{h,n}$, is defined for integers $h\ge 2$ and even numbers $n\ge 4$ as follows.

$$V(D_{h,n})=V(P_h\Box C_n),$$
$$E(D_{h,n})=E(P_h\Box C_n)\cup\{(i,j),(i,k)\}:i\in\{1,h\},\,j-k\ {\rm is\ odd}\}.$$
That is, $D_{h,n}$ is obtained from $P_h\Box C_n$ by extending the even cycles on the bottom and top layers (induced by $\{1\}\times V(C_n)$ and $\{h\}\times V(C_n)$) into complete bipartite graphs $K_{n/2,n/2}$.

For $i\in[h]$ and $j\in[n/2]$, we call the vertices $(i,j)$ and $(h+1-i,j+n/2)$ of $D_{h,n}$ a pair of \emph{opposite vertices}. The \emph{reduced drum} $D'_{h,n}$ is obtained from the drum $D_{h,n}$ by identifying all opposite pairs of vertices.
\end{defi}

Observe that every vertex of the drum $D_{h,n}$ has a unique opposite vertex and the involution $\gamma$ switching each opposite pair is an automorphism of $D_{h,n}$.

As opposite vertices are not adjacent, the map $\iota$ that brings a vertex
$v$ of $D_{h,n}$ to the vertex of $D'_{h,n}$ representing the pair of opposite
vertices containing $v$ is a graph homomorphism. For $h>2$ and $n>4$ opposite vertices do not have a common neighbor. In this case, as the switching $\gamma$ of opposite pairs is an automorphism of $D_{h,n}$, the homomorphism $\iota$ is, in fact, a double cover. We note that opposite vertices have no common neighbors in either $D_{2,n}$ if $4$ divides $n$ or in $D_{h,4}$ if $h$ is even. So the map $\iota$ is a double cover in these cases, too.

Here we constructed the reduced drum $D'_{h,n}$ as a factor of the drum $D_{h,n}$. It may be instructive to also construct it from ``half of the drum $D_{h,n}$'', that is from $P_{\lceil h/2\rceil}\Box C_n$, where the bottom layer is extended to a complete bipartite graph. In order to do so we only have to adjust opposite vertices in the top layer. For even $h$ we connect them, for odd $h$ we identify them as follows.

\begin{lemma}\label{ize}
For even integers $h\ge2$ and $n\ge4$ the reduced drum $D'_{h,n}$ can be obtained from $P_{h/2}\Box C_n$ by adding the edges $\{(1,i)(1,j):i-j\hbox{\rm\ is odd}\}\cup\{(h/2,i)(h/2,j):|i-j|=n/2\}$.

For odd $h\ge3$ and even $n\ge4$ the reduced drum $D'_{h,n}$ can be obtained from $P_{(h+1)/2}\Box C_n$ by adding the edges $\{(1,i)(1,j):i-j\hbox{\rm\ is odd}\}$ and then identifying the vertex $((h+1)/2,i)$ with the vertex $((h+1)/2,i+n/2)$ for all $1\le i\le n/2$.
\end{lemma}

\proof
We map the vertex $(i,j)$ with $i\le h/2$ of the graph constructed in the lemma to the vertex in $D'_{h,n}$ obtained from identifying the vertex $(i,j)$ of $D_{h,n}$ with its opposite vertex. In case $h$ is odd we further have the vertices $w_i$ obtained by identifying $((h+1)/2,i)$ with $((h+1)/2,i+n/2)$ for $1\le i\le n/2$. As $((h+1)/2,i)$ and $((h+1)/2,i+n/2)$ are opposite in $D_{h,n}$ they are also identified in $D'_{h,n}$ and the resulting vertex is where we map $w_i$. It is straightforward to check that the mapping so defined is an isomorphism. \qed
\smallskip

Looking at the reduced drum $D'_{h,n}$ as constructed in Lemma~\ref{ize} we see $\lceil h/2\rceil$ layers of the product $P_{\lceil h/2\rceil}\Box C_n$. The bottom and top layers (corresponding to the vertex $1$ and $\lceil h/2\rceil$ of $P_{\lceil h/2\rceil}$, respectively) are modified, the rest of the layers are copies of $C_n$. The bottom layer is a complete bipartite graph $K_{n/2,n/2}$. For odd $h$ the top layer is the cycle $C_{n/2}$ (or $K_2$ in case $n=4$), while for even $h$ it is the cycle $C_n$ with the opposite vertices connected. Using notation from \cite{GH}, we call this last graph a \emph{M\"obius ladder} and denote it by $M_n$. Only neighboring layers are connected and they are connected by a perfect matching connecting corresponding vertices except in case of the top layer for odd $h$, where each vertex is connected to two opposite vertices of the previous layer. The case $h=2$ is degenerate because the whole graph is a single layer that is both the bottom and the top layer. It is the complete bipartite graph $K_{n/2,n/2}$ if $n$ is not divisible by $4$, and in case $n$ is divisible by $4$ one has to add a perfect matching on both sides of $K_{n/2,n/2}$ to obtain $D'_{2,n}$.

Now we are ready to describe the structure of $4$-chromatic Schrijver
graphs. As mentioned in the Introduction, this description is
already implicit in the proof of the main result of Braun's paper
\cite{BBraun1}. It
is also given more explicitely (though in somewhat different terms)
in \cite{BBraun2} as well as in an
earlier unpublished manuscript by Li \cite{Li}.

\begin{thm}\label{thm:SG4struct} (cf. {\rm \cite{BBraun1, BBraun2, Li}})
For $k\ge 1$ the Schrijver graph ${\rm SG}(2k+2,k)$ is isomorphic to the reduced drum $D'_{k+1,2k+2}$.
\end{thm}

\begin{figure}[!htbp]
\centering
\includegraphics[scale=0.36]{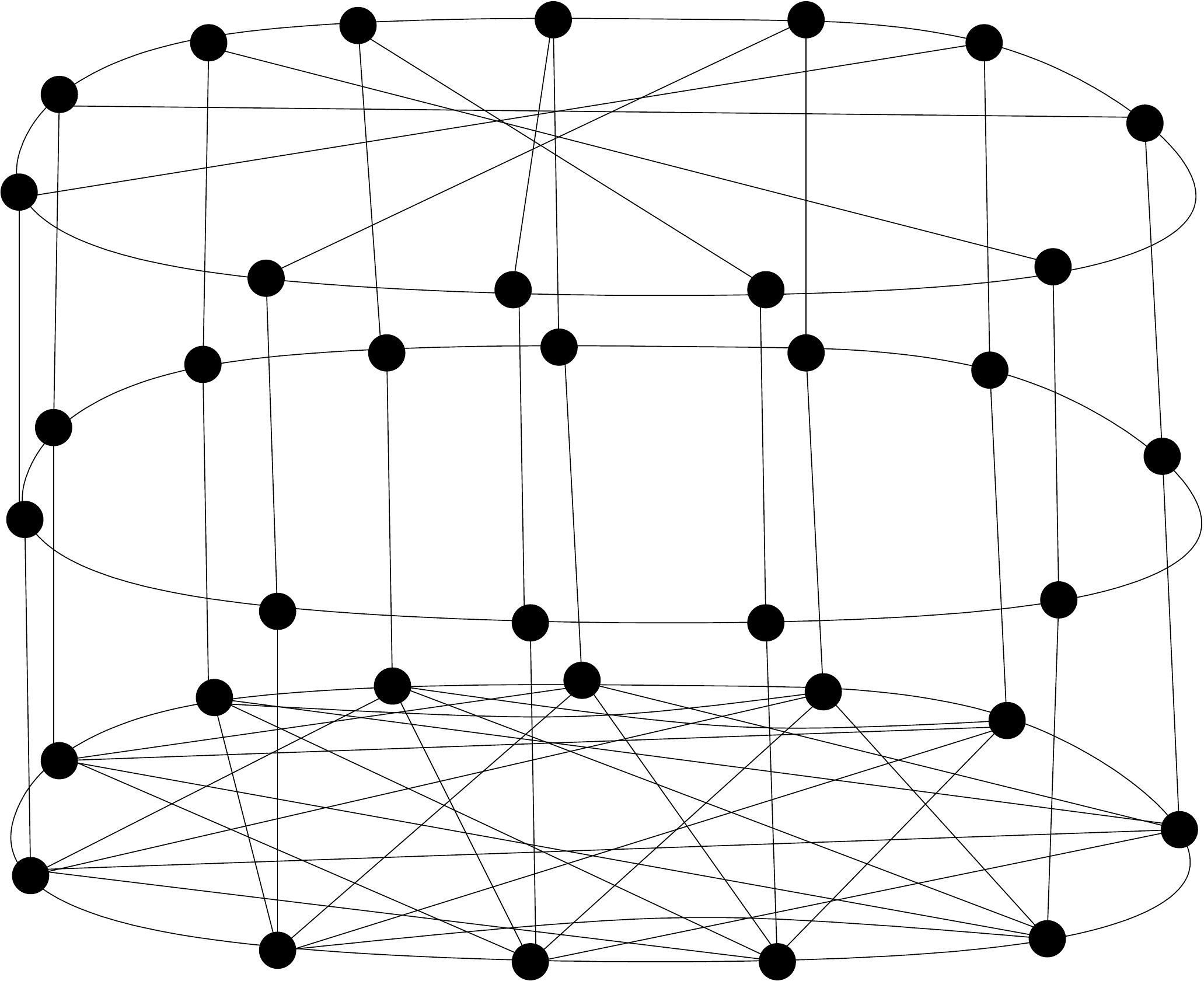}
\caption{This figure shows the graph ${\rm SG}(12,5)\cong D'_{6,12}$ the way
  we think about
  it: the ``bottom cycle'' induces a $K_{6,6}$, while the ``top cycle''
  induces the M\"obius ladder
$M_{12}$. See also Figures~\ref{fig:62}~and~\ref{fig:83} for more
symmetric pictures of the graphs ${\rm SG}(6,2)$ and ${\rm SG}(8,3)$,
respectively.}
\label{fig:125}
\end{figure}

\proof In case $k=1$ the theorem claims that $\SG(4,1)$ is isomorphic to $D'_{2,4}$. This holds as both of these graphs are complete graphs on four vertices. For the rest of the proof, we assume $k\ge2$ and set $C=C_{2k+2}$, $\SG=\SG(2k+2,k)$ and $D=D_{k+1,2k+2}$. All our calculations will be modulo $2k+2$. We denote the edge of $C$ connecting $a$ and $a+1$ by $e_a$.

Recall that the vertices of $\SG$ are the size $k$ independent sets in $C$. As $C$ is a cycle, each of the $k$ elements of a set $v\in V(\SG)$ is incident to two edges of $C$ and since $v$ is an independent set these $2k$ edges must all be distinct. Therefore, there are exactly two edges in $C\setminus v$. Let us define $e(v)=E(C\setminus v)$, so $e(v)=\{e_a,e_b\}$ for two distinct edges $e_a$ and $e_b$ of $C$. It is easy to see that for indices $a$ and $b$ of distinct parity there is exactly one independent set of size $k$ among the vertices of $C$ outside the endpoints of $e_a$ and $e_b$, while for distinct indices $a$ and $b$ of the same parity, there are none. Thus, $e(v)$ determines $v\in V(\SG)$ and $e(v)$ must be of the form $\{e_a,e_b\}$ with $a-b$ odd.

Consider the map $f:V(D)\to V(\SG)$ defined by $e(f((i,j)))=\{e_{i+j-1},e_{j-i}\}$. The function is well defined as $i+j-1$ and $j-i$ are of different parity.

We claim that $f$ is a cover. Consider an arbitrary vertex $(i,j)\in V(D)$. The neighbors of the vertex $v=f((i,j))$ in $\SG$ are the independent sets of size $k$ in $C\setminus v$. Here $C\setminus v$ has $k+2$ vertices and its edge set is $e(v)=\{e_{i+j-1},e_{j-i}\}$. Notice that these are adjacent edges if $i=1$ or $i=k+1$, but they are independent otherwise. Let us first assume $1<i<k+1$. Then $v$ has exactly $4$ neighbors, each obtained from $X=V(C)\setminus v$ by removing one endpoint of both edges in $e(v)$. Notice that these are exactly the vertices $f(w)$, where $w$ runs over the four neighbors of $(i,j)$. We have:
\begin{eqnarray*}
f((i-1,j))&=&X\setminus\{i+j-1,j-i+1\},\\
f((i+1,j))&=&X\setminus\{i+j,j-i\},\\
f((i,j-1))&=&X\setminus\{i+j-1,j-i\},\\
f((i,j+1))&=&X\setminus\{i+j,j-i+1\}.
\end{eqnarray*}
It remains to deal with the vertices $(i,j)$ of $D$ with $i\in\{1,k+1\}$. We may assume by symmetry that $i=1$. Now the vertex $v=f((1,j))$ has $e(v)=\{e_{j-1},e_j\}$ . Here $j$ is the common vertex of the two edges. Therefore $v$ has a single neighbor containing $j$: it is obtained by removing the two neighbors of $j$ in $C$ from $X=V(C)\setminus v$ and is the image of the only neighbor of $(1,j)$ outside its layer:
\begin{eqnarray*}
f((2,j))&=&X\setminus\{j-1,j+1\}.
\end{eqnarray*}
The rest of the neighbors of $v$ are
$k$-subsets of $S=X\setminus\{j\}$. Note that we have $S=\{u\in[2k+2]:u-j\hbox{ is odd}\}$. Any $k$-subset of $S$ is of the form $w=S\setminus\{j'\}$ for some $j'\in S$, therefore we have $e(w)=\{e_{j'-1},e_{j'}\}$, which means $w=f((1,j'))$ for the neighbor $(1,j')$ of $(1,j)$. This finishes the proof of the claim that $f$ is a covering.

It remains to find the inverse image of a vertex of $v$ in $\SG$. Let $e(v)=\{e_a,e_b\}$. We have $f((i,j))=v$ if and only if $i+j-1=a$ and $j-i=b$ or the other way around, $i+j-1=b$ and $j-i=a$. Recall that we still do calculations modulo $2k+2$, so these are systems of congruences. Both of these systems have two solutions modulo $2k+2$. Namely, if $i_0,j_0$ represent a solution, then the other three solutions will be $(k+2-i_0,j_0+k+1)$, $(i_0+k+1,j_0+k+1)$ and $(2k+3-i_0,j_0)$. A pair $(i,j)$ represent a vertex of $D$ only if $i\in[k+1]$. Thus, always two of the four solutions represent vertices in $D$, either the first two or the last two. Notice that they represent a pair of opposite vertices in $D$. This means that $f$ is a double cover that collapses the opposite pairs of vertices. To finish the proof of the theorem notice that the image of a cover $f$ is always isomorphic to the graph obtained by identifying the sets of vertices in the domain of $f$ that are mapped to the same vertex. We obtained $D'_{k+1,2k+2}$ from $D$ exactly like this: we identified the pairs of opposite vertices. \qed

\begin{figure}[!htbp]
\centering
\includegraphics[scale=0.50]{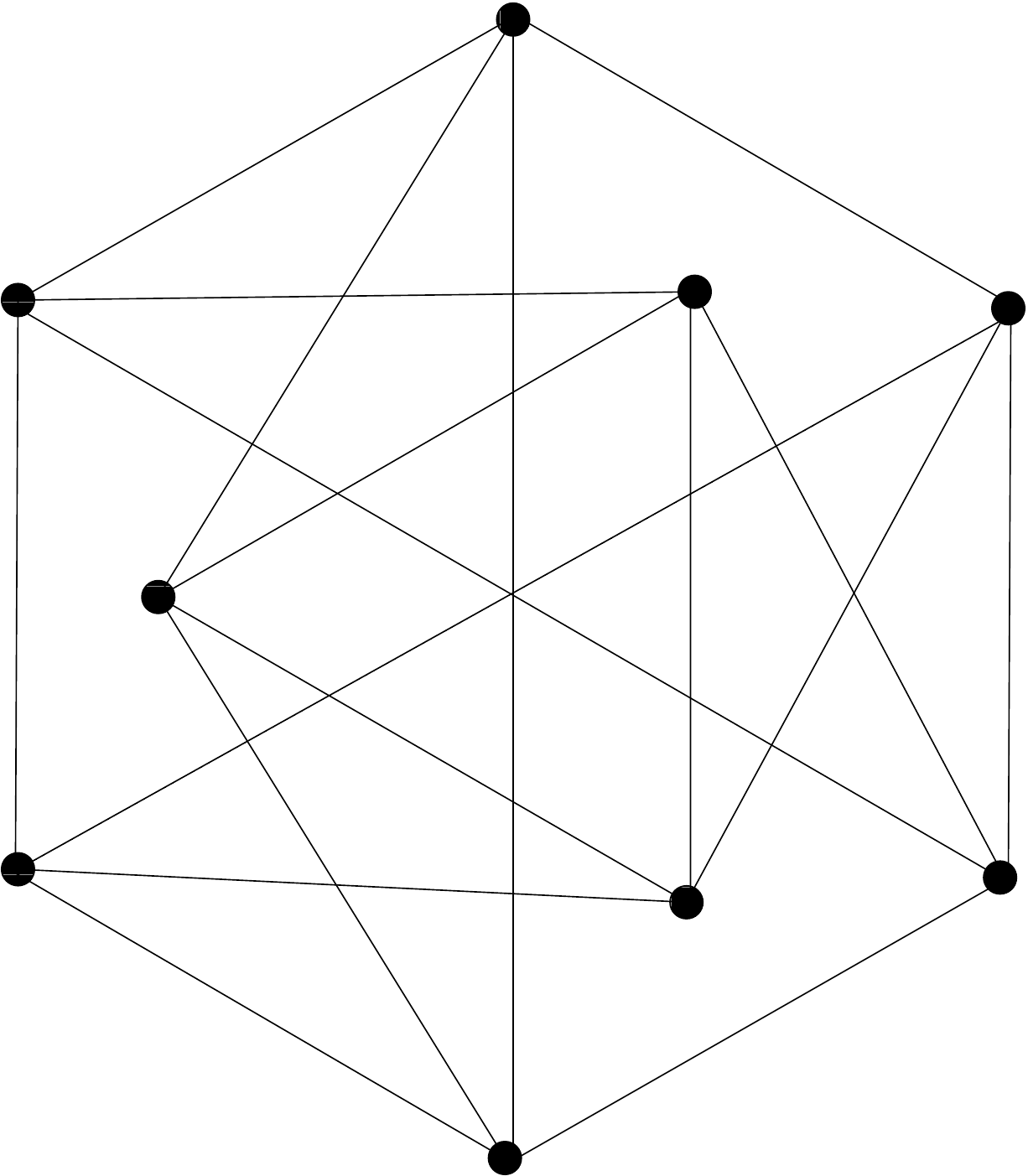}
\caption{The graph ${\rm SG}(6,2)$. This is the first ``interesting''
  Schrijver graph, that is, the smallest one which is not just a complete
  graph or an odd cycle. The ``inner cycle'' here is the ``top
  cycle'' $C_3$ obtained after identifying opposite points of a $C_6$, while the
``outer cycle'' is the $C_6$ extended to a complete bipartite graph $K_{3,3}$.}
\label{fig:62}
\end{figure}

\begin{figure}[!htbp]
\centering
\includegraphics[scale=0.50]{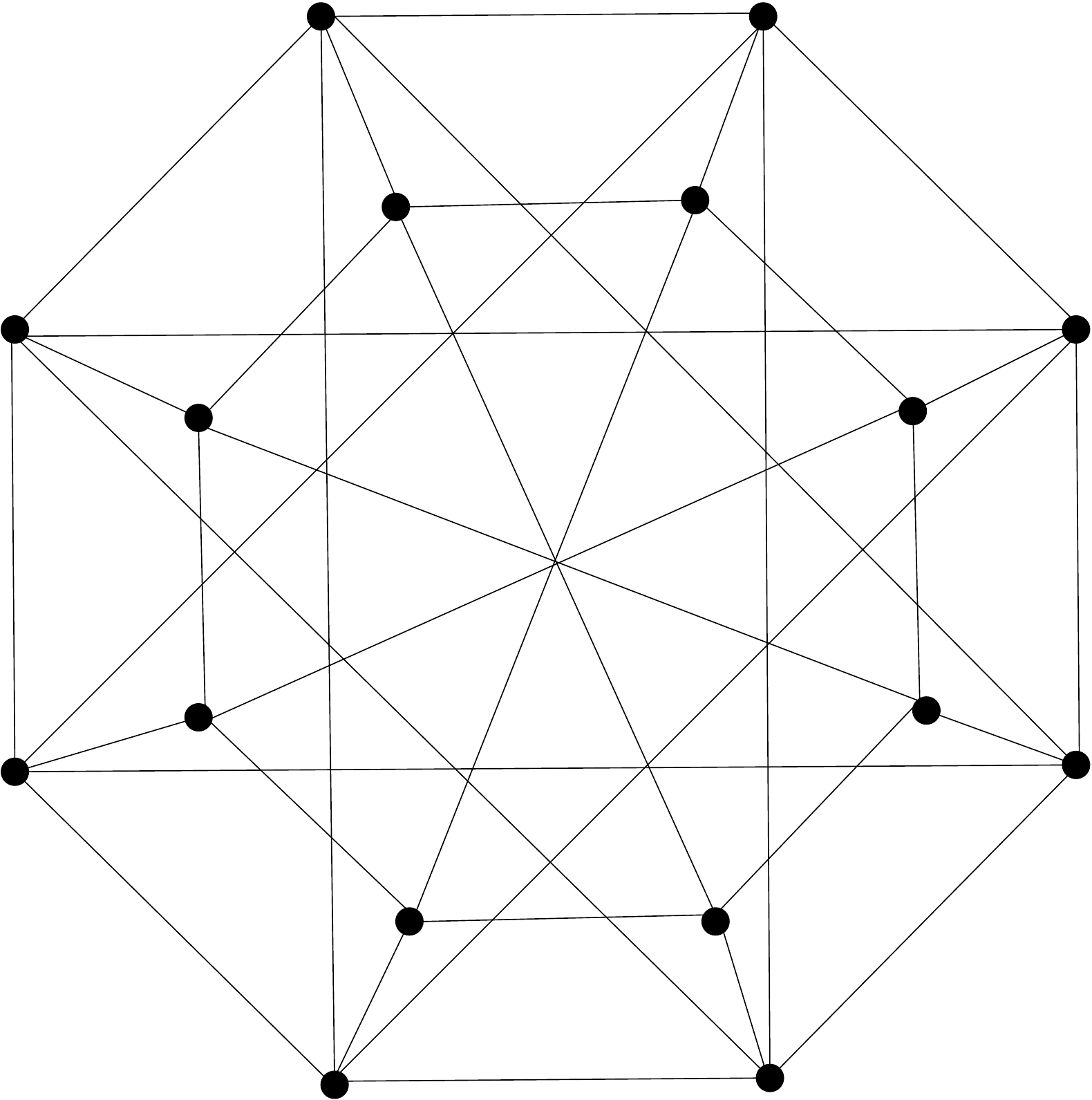}
\caption{The graph ${\rm SG}(8,3)$. We have drawn
  the complete bipartite graph $K_{4,4}$ as the ``outer cycle'', while the
  ``inner cycle'' induces the M\"obius ladder $M_8$.}
\label{fig:83}
\end{figure}

\section{The chromatic number of reduced drums}\label{sect:proof}

Our goal in this chapter is to prove the following statement in an elementary
way.

\begin{thm}\label{thm:reduceddrum}
For $h\ge2$ integer and $n\ge4$ even, we have $\chi(D'_{h,n})=2$ if $h+n/2$ is odd and $\chi(D'_{h,n})=4$ if $h+n/2$ is even.
\end{thm}
\medskip

\noindent
Note that by Theorem~\ref{thm:SG4struct} this implies that $\SG(2k+2,k)=4$ for all $k\ge1$.

The main part of Theorem~\ref{thm:reduceddrum} is the non-three-colorability
in case $h+n/2$ is even. We state that separately in Lemma~\ref{lem:non3col}
and give the simple proof for the rest of Theorem~\ref{thm:reduceddrum} right
after stating this lemma.

\begin{lemma}\label{lem:non3col}
There is no proper 3-coloring of the reduced drum $D'_{h,n}$ if $h\ge2$, $n\ge4$ even and $h+n/2$ is even.
\end{lemma}

\noindent\emph{Proof of Theorem~\ref{thm:reduceddrum} using Lemma~\ref{lem:non3col}.}
$D_{h,n}$ is always bipartite, with the color of the vertex $(i,j)$ in the proper 2-coloring given by the parity of $i+j$. This coloring gives identical colors to opposite vertices $(i,j)$ and $(h+1-i,j+n/2)$ whenever $h+n/2$ is odd. So in this case the factor $D'_{h,n}$ inherits the same 2-coloring.

If $h+n/2$ is even, then Lemma~\ref{lem:non3col} gives $\chi(D'_{h,n})\ge4$, so it is enough to find a proper $4$-coloring. We consider $D'_{h,n}$ as constructed in Lemma~\ref{ize}. We start with an arbitrary proper coloring of the top layer with colors from $\{1,2,3,4\}$. The top layer is a cycle $C_{n/2}$ or a M\"obius ladder $M_n$, so even three colors suffice in all cases except for $M_4\simeq K_4$. For vertices $(i,j)$ in all other layers, that is for $i<\lceil h/2\rceil$ we use a color of the same parity as $i+j$. For the vertices $(\lceil h/2\rceil-1,j)$ we choose from the two colors allowed by this parity rule to avoid the color of the only neighbor of this vertex in the top layer. For all other vertices we can choose arbitrarily and we obtain a proper 4-coloring. \qed
\medskip

\noindent
Our key tool in the proof of Lemma~\ref{lem:non3col} is the following notion.

\medskip
\par\noindent
\begin{defi}\label{defi:wind}
When speaking of a 3-coloring of a graph we will always assume that the colors are $0$, $1$ and $2$. Let $f:V(G)\to \{0,1,2\}$ be a proper $3$-coloring of a graph $G$. Consider an arbitrary edge $e=ab$ of $G$. We define $W_f(\overrightarrow e)\in\{-1,1\}$ for an orientation $\overrightarrow e$ of $e$ from $a$ to $b$ such that $W_f(\overrightarrow e)\equiv f(b)-f(a)$ modulo $3$. Note that this is possible as $f$ is a proper coloring. Let $C$ be a cycle in $G$ and let $\overrightarrow C$ be a cyclic orientation of $C$. We set
$$W_f(\overrightarrow C):=\frac{1}{3}\sum_{\overrightarrow e\in E(\overrightarrow C)} W_f(\overrightarrow e)$$
and
$$W_f(C):=|W_f(\overrightarrow C)|.$$
We call $W_f(C)$ the \emph{winding number} of $f$ on $C$.
\end{defi}

\smallskip
\par\noindent
Note that $\sum_{\overrightarrow e\in E(\overrightarrow C)} W_f(\overrightarrow e)\equiv0$ modulo 3 by the definition, so the winding number is an integer. Informally, it is the number of
times we wind around the triangle $K_3$ on vertices $0,1,2$ as we follow how the
homomorphism $f$ maps $C$ to $K_3$. The opposite orientation $\overleftarrow
C$ of $C$ yields $W_f(\overleftarrow C)=-W_f(\overrightarrow C)$, so the
winding number is well defined.

\medskip
\par\noindent
The proof will be obtained as a consequence of the following lemmata.

\begin{lemma}\label{lem:nonzero1}
If $f$ is a proper $3$-coloring of a graph $G$ and $C$ is an $m$-cycle in $G$, then $W_f(C)\equiv m$ modulo 2. We have $W_f(C)=0$ if $C$ is a 4-cycle.
\end{lemma}

\proof
The first statement immediately follows by observing that $3W_f(\overrightarrow C)$ is the sum of $m$ values $W_f(\overrightarrow e)$ and each of these values is $-1$ or $1$. The winding number is at most $m/3$, so in case $m=4$ it has to be even and less than 2, so it is 0.
\hfill$\Box$

\medskip
\par\noindent
\begin{lemma}\label{lem:ketkor}
Let $f$ be a proper $3$-coloring of the graph $P_2\Box C_m$ and let $C^{(1)}$
and $C^{(2)}$ denote its two layers. We have
$$W_f(C^{(1)})=W_f(C^{(2)}).$$
\end{lemma}

\proof
For $j\in[2]$, consider the orientation $\overrightarrow C^{(j)}$ of $C^{(j)}$ going as $(j,1)\to(j,2)\to(j,3)\to\cdots\to(j,m)\to(j,1)$. For $i\in[m]$, consider the oriented 4-cycle $\overrightarrow D^{(i)}$: $(1,i)\to (1,i+1)\to (2,i+1)\to (2,i)\to (1,i)$. Here $m+1$ is meant to be $1$. We have
$$W_f(\overrightarrow C^{(1)})-W_f(\overrightarrow C^{(2)})=\sum_{i=1}^mW_f(\overrightarrow D^{(i)}).$$
Indeed, when applying the definition of the winding number to the right hand side, the terms $W_f(\overrightarrow e)$ for the ``vertical'' edges $e=(1,i)(2,i)$ cancel because both orientations of these edges come up once. The terms corresponding to ``horizontal'' edges remain but the same terms show up on both sides of the equation.

From Lemma~\ref{lem:nonzero1} the right hand side of the above equation is zero, so we have $W_f(\overrightarrow C^{(1)})=W_f(\overrightarrow C^{(2)})$ and also $W_f(C^{(1)})=W_f(C^{(2)})$. \qed

\begin{cor}\label{cor:sokkor}
Let $f$ be a proper $3$-coloring of the graph $P_h\Box C_n$ and denote by
$C^{(j)}$ the layer induced by the vertices $\{(j,i): 1\le i\le n\}$.
Then $$W_f(C^{(1)})=W_f(C^{(h)}).$$
\end{cor}

\proof
Lemma~\ref{lem:ketkor} implies $W_f(C^{(j)})=W_f(C^{(j+1)})$ for
$j=1,\dots,h-1$, thus we have  the equalities
$W_f(C^{(1)})=W_f(C^{(2)})=\cdots =W_f(C^{(h)})$.
\hfill$\Box$

\medskip
\par\noindent
\begin{lemma}\label{lem:zero}
If $f$ is a proper $3$-coloring of a complete bipartite graph $K_{m,m}$, and
$C$ is any cycle of this complete bipartite graph, then
$W_f(C)=0$.
\end{lemma}

\proof
The statement immediately follows from the fact, that one full side of the complete
bipartite graph, say $H$, must be colored with one color if only $3$ colors are
used. Let us cut the cycle $C$ into $2$-edge paths at the vertices in $H$. At any cyclic orientation $\overrightarrow C$ of the cycle $C$ the contributions of the edges in one of these paths to $W_f(\overrightarrow C)$ cancel out, so we must have $W_f(\overrightarrow C)=0$ and therefore $W_f(C)$ must also be zero. \hfill$\Box$

\medskip
\par\noindent

\medskip
\par\noindent
\begin{lemma}\label{lem:nonzero2}
Let $n\ge4$ be even and $f$ a proper $3$-coloring of the cycle
$C_n$ that also puts different colors on vertices that are opposite to each other. In other words, $f$ remains a proper coloring even if we extend the cycle into a M\"obius ladder $M_n$ by adding the edges connecting the $n/2$ pairs of opposite vertices. Then $W_f(C_n)\equiv n+2$ modulo $4$.
\end{lemma}

\proof
As usual, we assume $V(C_n)=[n]$. When referring to the vertex $n+1$ we mean the vertex $1$. For $i\in[n]$, let $\overrightarrow e_i$ stand for the oriented edge from $i$ to $i+1$. Let $\overrightarrow C_n$ stand for the orientation of $C_n$ so given.

Consider the sums $m_j:=\sum_{i=j}^{j+n/2-1}W_f(\overrightarrow e_i)$ for $j\in[n/2+1]$. Clearly, $m_j\equiv f(j+n/2)-f(j)$ modulo 3. As $j$ and $j+n/2$ are opposite vertices in $C_n$ we have $f(j+n/2)\ne f(j)$ so $m_j$ is not divisible by 3. Also, as $m_j$ is the sum of $n/2$ terms, each $-1$ or $1$, we have $m_j\equiv n/2$ modulo 2. So if $n/2$ is even, then $m_j$ is congruent with $2$ or $4$ modulo $6$, while if $n/2$ is odd, then $m_j$ is congruent to $-1$ or $1$ modulo $6$. We can write $m_j=6s_j+t_j$, where $s_j$ is an integer and $t_j\in\{2,4\}$ for all $j$ if $n/2$ is even and $t_j\in\{-1,1\}$ for all $j$ if $n/2$ is odd.

For $j\in[n/2]$ we have $m_{j+1}-m_j=W_f(\overrightarrow e_{j+n/2})-W_f(\overrightarrow e_j)$ and therefore $|m_{j+1}-m_j|\le2$. This implies $s_{j+1}=s_j$. As this holds for all $j\in[n/2]$ we also have $s_1=s_{n/2+1}$.

By the definition of the winding number we have $3W_f(\overrightarrow C_n)=m_1+m_{n/2+1}$. Here $m_1+m_{n/2+1}=(6s_1+t_1)+(6s_{n/2+1}+t_{n/2+1})=12s_1+t_1+t_{n/2+1}$. We see that $t_1+t_{n/2+1}$ must be divisible by $3$. In case $n/2$ is even we have $t_1,t_{n/2+1}\in\{2,4\}$, so $t_1+t_{n/2+1}=6$ must hold. This yields $3W_f(\overrightarrow C_n)=12s_1+6$, so $W_f(\overrightarrow C_n)=4s_1+2$ as claimed. If $n/2$ is odd, we have $t_1,t_{n/2+1}\in\{-1,1\}$, so we must have $t_1+t_{n/2+1}=0$ and $W_f(\overrightarrow C_n)=4s_1$.
\hfill$\Box$

\medskip
\par\noindent
{\em Proof of Lemma~\ref{lem:non3col}.}
Assume that for some integer $h\ge2$ and even $n\ge4$ the reduced drum $D'_{h,n}$ admits a proper 3-coloring $f$. We need to prove that $h+n/2$ is odd.

Let us first assume that $h$ is even. By Lemma~\ref{ize} we can consider $D'_{h,n}$ as the graph $P_{h/2}\Box C_n$  with a few additional edges. Let $C^{\mathrm{bottom}}$ stand for the bottom layer of $P_{h/2}\Box C_n$ and $C^{\mathrm{top}}$ for its top layer.
The additional edges make the subgraph of $D'_{h,n}$ induced by $C^{\mathrm{bottom}}$ a complete bipartite graph, so by Lemma~\ref{lem:zero} we must have $W_f(C^{\mathrm{bottom}})=0$. By Corollary~\ref{cor:sokkor} we have $W_f(C^{\mathrm{top}})=W_f(C^{\mathrm{bottom}})$. The additional edges make the subgraph of $D'_{h,n}$ induced by $C^{\mathrm{top}}$ a M\"obius ladder $M_n$, so by Lemma~\ref{lem:nonzero2} we have $0=W_f(C^{\mathrm{top}})\equiv n+2$ modulo $4$. This means $n\equiv 2$ modulo 4, so $h+n/2$ is odd.

Now we assume that $h$ is odd. By Lemma~\ref{ize} we can consider $D'_{h,n}$
as constructed from $P_{(h+1)/2}\Box C_n$ in two steps. We refer to the bottom
and top layers of $P_{(h+1)/2}\Box C_n$ as $C^{\mathrm{bottom}}$ and
$C^{\mathrm{top}}$. In the first step we construct $\hat D$ from
$P_{(h+1)/2}\Box C_n$ by adding edges to $C^{\mathrm{bottom}}$ to make it into
a complete bipartite graph. Then we identify the opposite vertices in
$C^{\mathrm{top}}$ to obtain $D'_{h,n}$. Clearly, the 3-coloring $f$ can be
extended to a proper 3-coloring $\hat f$ of $\hat D$ by keeping the color of the
vertices outside $C^{\mathrm{top}}$ and assigning both opposite vertices in
$C^{\mathrm{top}}$ the color that $f$ assigns to the new vertex resulting from
their identification. As in the previous case we must have $W_{\hat
  f}(C^{\mathrm{top}})=W_{\hat f}(C^{\mathrm{bottom}})=0$ by
Corollary~\ref{cor:sokkor} and Lemma~\ref{lem:zero}. Let $C^*$ stand for the cycle of length $n/2$
obtained from the cycle $C^{\mathrm{top}}$ by identifying the opposite vertices. In case $n=4$ the graph $C^*$ is not a cycle but then $h+n/2$ is odd anyway. We clearly have $W_{\hat f}(C^{\mathrm{top}})=2W_f(C^*)$ as each contribution $W_f(\overrightarrow e)$ to the winding number of an orientation of $C^*$ appears twice (on two opposite edges) as the contribution to the winding number of the corresponding orientation of $C^{\mathrm{top}}$. So we must have $W_f(C^*)=0$, and by Lemma~\ref{lem:nonzero1} this means that $n/2$ is even. This makes $h+n/2$ odd again and completes the proof. \qed

Finally we state a similar result that is not directly related to Schrijver
graphs but follows easily from the lemmata above. For $h\ge2$ integer and even
$n\ge4$ let us construct the graph $X_{h.n}$ from $P_h\Box C_n$ by adding the
edges $(1,i)(1,i+n/2)$ for $i\in[n/2]$ and identifying the vertices $(h,i)$
and $(h,i+n/2)$ for all $i\in[n/2]$. Note that $X_{h,n}$ is a $4$-regular graph and in case $n/2$ is odd, it is a subgraph of $D'_{2h-1,n}$.

\begin{thm}\label{bigyo}
For all $h\ge2$ and even $n\ge4$ we have $\chi(X_{h,n})=4$.
\end{thm}

\proof We write $C^{\mathrm{bottom}}$ and $C^{\mathrm{top}}$ for the bottom and top layers of $P_h\Box C_n$ and, in case $n\ge6$, let $C^*$ be the $n/2$-cycle of $X_{h,n}$ obtained by identifying the opposite vertices in $C^{\mathrm{top}}$. Let $\hat X$ be $P_h\Box C_n$ plus the additional edges making $C^{\mathrm{bottom}}$ into the M\"obius ladder $M_n$.

Assume for a contradiction that $f$ is a proper 3-coloring of $X_{h,n}$. We can extend $f$ to $\hat f$, a proper $3$-coloring of $\hat X$ that assigns the same color to opposite vertices of $C^{\mathrm{top}}$ as $f$ assigned to the vertex obtained by identifying them. We have $W_{\hat f}(C^{\mathrm{bottom}})\equiv n+2$ modulo 4 by Lemma~\ref{lem:nonzero2}. We have $W_{\hat f}(C^{\mathrm{top}})\equiv n$ modulo 4 if $n\ge 6$ because $W_{\hat f}(C^{\mathrm{top}})=2W_f(C^*)$ and $W_f(C^*)\equiv n/2$ modulo 2 by Lemma~\ref{lem:nonzero1}. In case $n=4$ Lemma~\ref{lem:nonzero1} gives $W_{\hat f}(C^{\mathrm{top}})=0$, so the same congruence $W_{\hat f}(C^{\mathrm{top}})\equiv n$ modulo 4 also holds. Therefore we have $W_{\hat f}(C^{\mathrm{top}})\ne W_{\hat f}(C^{\mathrm{bottom}})$ contradicting Corollary~\ref{cor:sokkor}. The contradiction proves $\chi(X_{h,n})\ge4$.

As we noted, $X_{h,n}$ is a subgraph of $D'_{2h-1,n}$ if $n/2$ is odd. In this
case $\chi(X_{h,n})\le\chi(D'_{2h-1,n})=4$ by
Theorem~\ref{thm:reduceddrum}. It remains to prove $\chi(X_{h,n})\le4$ if
$n/2$ is even.
In that case ignoring the edges added to $C^{\mathrm{bottom}}$ to make it a
M\"obius ladder the graph is bipartite. These extra edges form a matching, so another bipartite graph. Therefore $X_{h,n}$, as the union of two bipartite graphs, has chromatic number at most $4$.
\qed

\section{Critical edges} \label{sect:crit}

In this section we study when will the deletion of a single edge of a
Schrijver graph reduce its chromatic number. Here we consider arbitrary
Schrijver graphs, not only 4-chromatic ones, but point out how our results
apply to 4-chromatic Schrijver graphs. (To emphasize that the vertices of
Schrijver graphs are subsets of $[n]$, in what follows they will be denoted by capital
letters.)

\begin{defi}
We call the edge $e$ of a graph $G$ \emph{color-critical} or simply \emph{critical} if the removal of that edge decreases the chromatic number, that is $\chi(G\setminus\{e\})<\chi(G)$.

We call an edge $XY\in E({\rm SG}(n,k))$ \emph{interlacing} if the elements of $X$ and $Y$ alternate on the cycle $C_n$.
\end{defi}

We will see that distinguishing interlacing and non-interlacing edges of the Schrijver graphs is {\em not} enough to characterize critical edges, but it nevertheless seems to be an important distinction. In Section~\ref{otegy} we will show that non-interlacing edges are not critical, in Section~\ref{otketto} we show that several categories of interlacing edges are critical. In fact, we will see in Section~\ref{otharom} that our results cover all edges in the 4-chromatic Schrijver graphs that are the main focus of this paper. In Section~\ref{otnegy} we characterize the critical edges in the Schrijver graphs ${\rm SG}(n,2)$ and will see that for $n\ge8$ some of the interlacing edges are not critical.

\subsection{Critical edges are interlacing}\label{otegy}

As a preparation for our proof that all critical edges are interlacing, we
recall a version of B\'ar\'any's proof \cite{Barany} of the
Lov\'asz-Kneser theorem and also the proof of Schrijver's result
\cite{Sch78}. We start with a definition of the following infinite graphs
first defined by Erd\H{o}s and Hajnal \cite{EH}.

\begin{defi}
Let $S^d$ stand for the unit sphere around the origin in $\reals^{d+1}$. We
consider $\reals^{d+1}$ (hence also $S^d$) equipped with the Euclidean
distance, so the diameter of $S^d$ is $2$. The open \emph{hemisphere} around a
point $x\in S^d$ is the set
$\{y\in S^d\mid x\cdot y>0\}=$
$\{y\in S^d\mid|x-y|<\sqrt2\}$. (Here $x\cdot y$ is meant to be the scalar
product of the vectors from the origin to $x$ and $y$, respectively.) The \emph{Borsuk graph} $B(d,\varepsilon)$ is the infinite graph whose vertices are the points of $S^d$ and two vertices are connected by an edge if $|x+y|<\varepsilon$ (that is $y$ is $\varepsilon$-close to $-x$). Here $\varepsilon$ is a positive real.
\end{defi}

The following proposition is an equivalent form of the Borsuk-Ulam theorem
(cf. Lov\'asz \cite{LLgomb}):

\begin{prop}\label{borsukgraph}
For any $d\ge0$ integer and $\varepsilon>0$ real we have $\chi(B(d,\varepsilon))\ge d+2$.
\end{prop}

B\'ar\'any's proof uses the following lemma of Gale \cite{Gal}:

\begin{lemma}\label{gale}
For any integers $d,k>0$ there is a set of $d+2k$ points on $S^d$ such that each open hemisphere contains at least $k$ of these points.
\end{lemma}

Although B\'ar\'any's proof of the Lov\'asz-Kneser theorem used another form of the Borsuk-Ulam theorem and did not refer to the Borsuk graph, it can be reformulated as the next proposition. The statement ${\rm KG}(n,k)\ge n-2k+2$ then follows through Proposition~\ref{borsukgraph}.

\begin{prop}\label{BLK}
For positive integers $n\ge2k$ one can find $\varepsilon>0$ and a graph homomorphism from $B(n-2k,\varepsilon)$ to ${\rm KG}(n,k)$.
\end{prop}

\proof Using Lemma~\ref{gale} we find points $x_1,\dots,x_n\in S^d$ for $d=n-2k$ such that each open hemisphere contains at least $k$ of these points. For $x\in S^d$ let $f(x)$ stand for the $k$'th smallest distance $|x-x_i|$ for $1\le i\le n$. By the choice of the points $x_i$ we have $f(x)<\sqrt2$ for all $x\in S^d$. As $f$ is a continuous function on a compact domain it takes a maximum value $a<\sqrt2$. Now choose $\varepsilon=\sqrt2-a>0$. We map a vertex $x$ of $B(d,\varepsilon)$ to any vertex $V$ of ${\rm KG}(n,k)$ satisfying that for all elements $i\in V$ we have $|x-x_i|\le a$. Such a vertex $V$ exists since $f(x)\le a$. This map is a graph homomorphism from $B(d,\varepsilon)$ to ${\rm KG}(n,k)$ because if $xy$ is an edge of $B(d,\varepsilon)$ and we map $x$ to $V$ and $y$ to $W$, then $x_i$ for $i\in V$ is contained in the hemisphere around $x$, while $x_j$ for $j\in W$ is contained in the hemisphere around $-x$. These two hemispheres are disjoint so $V$ and $W$ must also be disjoint making $VW$ an edge in the Kneser graph. \qed

Schrijver \cite{Sch78} found the following stronger version of Gale's lemma:

\begin{lemma}\label{galesch}
For integers $d,k>0$ there are points $x_1,\ldots,x_{d+2k}\in S^d$ such that one can find a vertex $V$ in ${\rm SG}(d+2k,k)$ for each open hemisphere $H$ such that $x_i\in H$ for all $i\in V$.
\end{lemma}

Using this stronger lemma one can make the statement of Proposition~\ref{BLK} stronger as follows. Note that this implies Schrijver's result through Proposition~\ref{borsukgraph}. Although Schrijver did not refer to Borsuk graphs in his proof, this is the essence of his proof.

\begin{prop}\label{sch}
For positive integers $n\ge2k$ one can find $\varepsilon>0$ and a graph homomorphism from $B(n-2k,\varepsilon)$ to ${\rm SG}(n,k)$.
\end{prop}

\proof This proof is essentially identical to the proof of
Proposition~\ref{BLK}, the only difference is that we use Lemma~\ref{galesch} in
place of Lemma~\ref{gale}. We choose $x_1,\dots,x_n\in S^d$ for $d=n-2k$ to
satisfy the statement of Lemma~\ref{galesch}. For $x\in S^d$ we define $f(x)=\min_{V\in V(\mathrm{SG}(n,k))}(\max_{i\in V}|x-x_i|)$. As before $f(x)<\sqrt2$ by the choice of the points $x_i$ and since $f$ is continuous on a compact domain we have $f(x)\le a$ for some fixed $a<\sqrt2$ and all $x\in S^d$. We choose $\varepsilon=\sqrt2-a>0$ as before. We map a vertex $x$ of $B(d,\varepsilon)$ to any vertex $V$ of ${\rm SG}(n,k)$ with $|x-x_i|\le a$ for all $i\in V$. Just as in the proof of Proposition~\ref{BLK} such a vertex $V$ exists for all $x$ and the resulting map is a graph homomorphism. \qed

We can state an even stronger form of Gale's lemma as follows. Gale's lemma is
stated as Lemma~3.5.1 in Matou\v sek's book \cite{Matbook}. The proof there
(attributed to an observation by Ziegler) explicitly shows this stronger form.

\begin{lemma}\label{galeinter}
For integers $d,k>0$ there are points $x_1,\ldots,x_{d+2k}\in S^d$ such that for every $x\in S^d$ one can find an interlacing edge $VW$ in ${\rm SG}(d+2k,k)$ such that $x_i$ is contained in the open hemisphere around $x$ for every $i\in V$ while $x_j$ is contained in the open hemisphere around $-x$ for each $j\in W$.
\end{lemma}

At first glance one would think that with this strong form of Gale's lemma we can prove that the subgraph of a Schrijver graph formed by the interlacing edges has the same chromatic number as the entire Schrijver graph. Unfortunately, the lemma only ensures that in a homomorphism of a Borsuk graph to the Schrijver graph we can make the image of {\emph{any one edge} to be an interlacing edge, but we cannot ensure that the image of {\emph{all edges} are interlacing. Furthermore, the strong statement alluded to above is false in general. Removing \emph{all} the non-interlacing edges typically drastically reduces the chromatic number. The subgraph of ${\rm SG}(n,k)$ formed by the interlacing edges is called the \emph{interlacing graph} ${\rm IG}_{n,k}$. Litjens, Polak, Sevenster and Vena \cite{LPSV} show that its chromatic number is $\lceil n/k\rceil$.

Nevertheless, Lemma~\ref{galeinter} is enough to prove that none of the non-interlacing edges is critical in a Schrijver graph. The proof is through the following proposition.

\begin{prop}\label{interhom}
For positive integers $n\ge2k$ and a vertex $V_0$ of ${\rm SG}(n,k)$ one can find $\varepsilon>0$ and a graph homomorphism from $B(n-2k,\varepsilon)$ to the graph obtained from ${\rm SG}(n,k)$ by removing all non-interlacing edges incident to $V_0$.
\end{prop}

\proof Once again, we choose $x_1,\dots,x_n\in S^d$ for $d=n-2k$, this time to satisfy the statement of Lemma~\ref{galeinter}. For $x\in S^d$ we define $$f(x)=\min_{VW}(\max(\max_{i\in V}|x-x_i|,\max_{j\in W}|-x-x_j|)),$$ where the minimum is for interlacing edges $VW$ in ${\rm SG}(n,k)$. As before $f(x)<\sqrt2$ by the choice of the points $x_i$ and since $f$ is continuous on a compact domain we have $f(x)\le a$ for some fixed $a<\sqrt2$ and all $x\in S^d$.

We choose $\varepsilon=(\sqrt2-a)/2>0$ and set $b=a+\varepsilon=\sqrt2-\varepsilon$. We map a vertex $x$ of $B(d,\varepsilon)$ to a vertex $V$ of ${\rm SG}(n,k)$ such that $|x-x_i|\le b$ for all $i\in V$. Such a vertex exists by the choice of the points $x_i$. This time we do not choose an arbitrary such vertex $V$ but follow this order of priorities: We map $x$ to $V_0$ only if no other choice is available and we map $x$ to a non-interlacing neighbor of $V_0$ only if all other available choices are also non-interlacing neighbors of $V_0$ or possibly $V_0$ itself.

To see that this defines a graph homomorphism, we need to show that if $x$ is
mapped to $V$ and $y$ is mapped to $W$ for an edge $xy$ of $B(d,\varepsilon)$,
then $VW$ is an edge of ${\rm SG}(n,k)$, and further if $V=V_0$ or $W=V_0$,
then $VW$ is interlacing. As before, $x_i$ for $i\in V$ is contained in the
hemisphere around $x$, while $x_j$ for $j\in W$ is contained in the hemisphere
around $-x$, so $V$ and $W$ must be disjoint and thus an edge of ${\rm
  SG}(n,k)$. Assume now that $V=V_0$ (the case $W=V_0$ is symmetric). We need
to prove that $VW$ is an interlacing edge. By $f(x)\le a$ we have an
interlacing edge $V'W'$ in ${\rm SG}(n,k)$ with $|x-x_i|\le a$ for $i\in V$
and $|-x-x_j|\le a$ for $j\in W$. This makes $V'$ available as the image of
$x$ and $W'$ available as the image of $y$. (Note that we may have $y$ at
distance $\varepsilon$ away from $-x$, but even then the furthest point of
$W'$ from $y$ is still within distance $b$.) We map $x$ to $V_0$ only as a last resort, therefore we must have $V'=V_0$. But then $W'$ is an interlacing neighbor of $V_0$ where $y$ could be mapped, therefore $y$ cannot be mapped to a non-interlacing neighbor of $V_0$. This makes $VW$ an interlacing edge. \qed

\begin{thm}\label{noncritical}
No non-interlacing edge of a Schrijver graph is critical, further if we remove
all non-interlacing edges incident to a fixed vertex of ${\rm SG}(n,k)$ the chromatic number of the remaining graph is still $n-2k+2$.
\end{thm}

\proof This follows directly from Propositions~\ref{borsukgraph} and \ref{interhom}. \qed

Let us remark that with a more involved order of priorities among the vertices of the Schrijver graph one can extend the size of the set of non-interlacing edges that can be removed without effecting the chromatic number. But finding an edge-critical subgraph this way does not seem to be within reach.

\subsection{Many interlacing edges are critical}\label{otketto}

Here we show several classes of interlacing edges of Schrijver graphs that are critical. To show that an edge $VW$ of the Schrijver graph ${\rm SG}(n,k)$ is critical one has to give a proper $(n-2k+1)$-coloring of the graph ${\rm SG}(n,k)\setminus\{VW\}$, or equivalently an $(n-2k+1)$-coloring of the vertices of ${\rm SG}(n,k)$ in which only the edge $VW$ is monochromatic. For this we will start with the following coloring.

\begin{defi}
We call a (non-proper) $(n-2k+1)$-coloring of the Schrijver graph
${\rm SG}(n,k)$ a \emph{basic coloring} for the edge $VW$ if $n-2k$ of the
colors used are identified with the $n-2k$ elements of $[n]\setminus(V\cup W)$
and the color of a vertex $U$ is selected from $U\setminus(V\cup W)$ whenever
this set is not empty and the color of $U$ is the $(n-2k+1)$th color, called
\emph{special color} if $U\setminus(V\cup W)$ is empty.
\end{defi}

\begin{prop}\label{basic}
In a basic coloring of a Schrijver graph for an edge $VW$ we have that:
\begin{enumerate}
\item Both endpoints of all the monochromatic edges are colored the special color and
  these edges form and induced matching.
\item If $C_n[V\cup W]$ is connected or has two components, both on an odd
  number of vertices, then $VW$ is the only monochromatic edge in the coloring making it a critical edge in ${\rm SG}(n,k)$.
\end{enumerate}
\end{prop}

\proof Two vertices in the same non-special color class $x$ both contain $x$, so they are intersecting. Two vertices in the special color class form an edge if they partition $V\cup W$, proving part~1.

The vertices of ${\rm SG}(n,k)$ are independent $k$-sets in $C_n$. If
$C_n[V\cup W]$ has a single component (a path or the entire cycle $C_n$ in the
trivial case $n=2k$), then it has a single partition into independent sets,
while if it has two components (necessarily paths), then it has two such
partitions. If both of those two components contain an odd number of vertices,
then the parts in one of these partitions are unequal, so only one of the
partitions correspond to a monochromatic edge of ${\rm SG}(n,k)$, so $VW$ is the only monochromatic edge. This makes $VW$ critical and proves part~2 as the basic coloring uses $n-2k+1$ colors. \qed

\begin{defi}
We call a vertex $V$ of the Schrijver graph ${\rm SG}(n,k)$ \emph{regular} if $V$ contains at least one of any three consecutive vertices on the cycle $C_n$.
\end{defi}

Note that for $n>3k$ there are clearly no regular vertices in ${\rm SG}(n,k)$, but for a fixed chromatic number $n-2k+2=d+2$ and large $k$, the Schrijver graph ${\rm SG}(2k+d,k)$ has $\Theta(k^d)$ vertices and all but only $\Theta(k^{d-1})$ of them are regular.

\begin{prop}\label{regular}
 If $V$ is a regular vertex of the Schrijver graph ${\rm SG}(2k+d,k)$, then the degree of $V$ is $2^d$ and all the $2^d$ edges incident to $V$ are interlacing.
 \end{prop}

 \proof
 Let $n=2k+d$. Any vertex $V$ of ${\rm SG}(n,k)$ is an independent set of size
 $k$ in $C_n$, so $2k$ of the $n$ edges of $C_n$ is incident to an element of
 $V$ leaving $n-2k=d$ edges in $C_n\setminus V$. Note that we have already
 made this observation in the 4-chromatic case in the proof of
 Theorem~\ref{thm:SG4struct}. A neighbor of $V$ in ${\rm SG}(n,k)$ is an
 independent set of size $k$ in $C_n$ that is disjoint from $V$, that is, it is
 obtained by removing $d$ elements from the set $[n]\setminus V$ in such a way that we remove at least one endpoint of all the $d$ edges in $C_n\setminus V$.

Now assume that $V$ is regular. In this case the $d$ edges in $C_n\setminus V$
form a matching. To obtain a neighbor of $V$ we then must remove exactly one
endpoint of these $d$ edges from $[n]\setminus V$ in an arbitrary way, hence
the $2^d$ possibilities. Clearly, all these neighbors form interlacing edges
with $V$. \qed

 \begin{thm}\label{2regular}
 Any edge connecting two regular vertices in a Schrijver graph is critical.
 \end{thm}

\proof Let $VW$ be an edge in the Schrijver graph ${\rm SG}(n,k)$ connecting two regular vertices. By symmetry among the nodes of $C_n$ we may assume $1\in V$. Let $D=[n]\setminus(V\cup W)$. $D$ consist of $n-2k$ nodes of $C_n$. Among any three consecutive nodes on $C_n$ at least one is contained in $V$ and another one in $W$ (to make these vertices regular), so at most one of the three nodes is in $D$. This means that two elements of $D$ are separated by at least two nodes not in $D$ along the cycle $C_n$.

We start with a basic coloring $f_0$ of ${\rm SG}(n,k)$ for the edge $VW$ in which we choose the highest possible index whenever we have a choice. That is, we have $f_0(U)=\max(U\cap D)$ for vertices $U$ of ${\rm SG}(n,k)$ intersecting $D$ and $f_0(U)$ is the special color if $U$ and $D$ are disjoint.

We will obtain the coloring $f$ of the vertices of ${\rm SG}(n,k)$ by modifying $f_0$ at one end of each monochromatic edge except for $VW$ to make those other edges non-monochromatic. We will also make sure not to create new monochromatic edges in the process. By Proposition~\ref{basic}, all monochromatic edges form an induced matching, so the modifications in the coloring do not interact: it is enough if we find a non-special color that is not represented among the $f_0$-colors of the neighbors of one of the endvertices of the monochromatic edges considered.

Let therefore $XY$ be a monochromatic edge of $f_0$ other than $VW$. The sets $X$ and $Y$ partition $V\cup W$, so one of them, say $X$, contains $1$. As $XY$ is different from $VW$, neither of the sets $X\setminus V$ or $Y\setminus W$ is empty. Let $a=\min(X\setminus V)$ and $b=\min(Y\setminus W)$ and note that $a$ and $b$ are distinct and at least $2$. Let us consider the case $b<a$ first. In this case we recolor $X$ by setting $f(X)=b-1$. We need to show first that $b-1\in D$, so this is an available non-special color. The independent set $Y$ contains $b$, so it cannot contain $b-1$. We have $b\notin W$, so by $X\cup Y=V\cup W$ we have $b\in V$ and therefore $b-1\notin V$. This makes $b-1\notin X$, as otherwise we had $b-1\in X\setminus V$, making $a<b$, a contradiction. So $b-1\notin X\cup Y=V\cup W$ and therefore $b-1\in D$ as claimed.

It remains to show that $X$ does not have a neighbor $Z$ with $f_0(Z)=b-1$. Let us assume that such a neighbor $Z$ does exist. Let $D'=[n]\setminus(X\cup Z)$.

We have $b-1=f_0(Z)=\max(Z\cap D)$, so $Z$ avoids all nodes $d\in D$ with
$d>b-1$. These nodes are then contained in $D'$ as $X$ is disjoint from
$D$. Let us consider now $d\in D$ with $d<b-1$. We have $d\notin W$, so as $W$
is regular we must have a neighbor (in $C_n$) $d'\in W$ with $d'=d-1$ or $d'=d+1$. This makes $d,d'\notin V$ and as $a>b$ we must also have $d,d'\notin X$. As $Z$ is an independent set in $C_n$, it cannot contain both $d$ and $d'$, hence either $d$ or $d'$ shows up in $D'$. Finally consider $b-1\in D$. As $f_0(Z)=b-1$, $Z$ must contain $b-1$, so it contains neither of its neighbors in $C_n$. One of these neighbors is $b$ and $b\in Y$, implying $b\notin X$ and $b\in D'$. We have $b,b-1\notin W$ and $W$ is regular so it must contain $b-2$, the other neighbor of $b-1$ (note that $b-1\in D$, therefore $b-1>1$). From $b-2\in W$ we conclude $b-2\notin V$ and since $a>b$ we also have $b-2\notin X$ and thus $b-2\in D'$.

In summary, every element $d\in D$ is either contained in $D'$ or a neighbor of $d$ is contained in $D'$, and for element $b-1$ of $D$ both of its neighbors are contained in $D'$. We saw that elements of $D$ are separated by at least two other nodes along $C_n$, so there is no double counting, we must have $|D'|>|D|$. This is a contradiction, as we clearly have $|D'|=|D|=n-2k$. The contradiction proves that the recoloring of $X$ does not create new monochromatic edges.

We can deal with the case $a<b$ similarly: we recolor $Y$ by setting $f(Y)=a-1$. We can show that $a-1\in D$ and that this does not create a monochromatic edge the same way as we did in the case $a>b$. After dealing with each monochromatic edge of $f_0$ other than $VW$ (and defining $f(U)=f_0(U)$ for all other the vertices $U$) we obtain the coloring $f$ of the vertices of ${\rm SG}(n,k)$ with $n-2k+1$ colors, where $VW$ is the only monochromatic edge showing that $VW$ is critical.
\qed

\subsection{Critical edges of $\SG(2k+2,k)$}\label{otharom}

In this short section we show that in case of the 4-chromatic Schrijver graphs (that are the focus of this paper) the results in Sections~\ref{otegy} and \ref{otketto} are enough to characterize all edges with respect to color-criticality.

\begin{thm}\label{characterize}
An edge of the Schrijver graph $\SG(2k+2,k)$ is color-chritical if and only if it is interlacing.
\end{thm}

\proof The only if part of the theorem is true for all Schrijver graphs as stated in Theorem~\ref{noncritical}.

For the if part consider an interlacing edge $VW$ in $\SG(2k+2,k)$. The set
$D=[n]\setminus(V\cup W)$ has two elements. Therefore $C_n[V\cup
  W]=C_n\setminus D$ is either a path or the disjoint union of two paths. In
the former case or if there are two paths, both on an odd number of vertices,
then $VW$ is critical by Proposition~\ref{basic}. If, however, there are two
paths and both have at least $2$ vertices, then $V$ and $W$ are regular
(otherwise they cannot be interlacing), so $VW$ is critical by Theorem~\ref{2regular}. We have covered all the possible cases as $C_n[V\cup W]$ has $2k$ vertices, so if it consists of two components, then the parity of their sizes agree. \qed

Recall that the Schrijver graph $\SG(2k+2,k)$ is isomorphic to the reduced
drum $D'_{k+1,2k+2}$ by Theorem~\ref{thm:SG4struct} and the latter can be
obtained from $P_{\lfloor k/2\rfloor+1}\Box C_{2k+2}$ by slightly modifying
the top and bottom layers as described in Lemma~\ref{ize}. This gives another (more geometric) view of the 4-chromatic Schrijver graphs and it is instructive to identify the critical edges in this view of them.

\begin{prop}\label{geomch}
In the view of the reduced drum $D'_{k+1,2k+2}$ given in Lemma~\ref{ize}, the regular vertices of the isomorphic graph $\SG(2k+2,k)$ correspond to exactly the vertices outside the bottom layer.

The non-critical edges of the reduced drum $D'_{k+1,2k+2}$ are exactly the edges we added to the bottom layer to make it a complete bipartite graph. That is, all edges not contained in the bottom layer are critical and so are also the edges of the copy of the cycle $C_{2k+2}$ that is the bottom layer of the graph $P_{\lfloor k/2\rfloor+1}\Box C_{2k+2}$. All remaining edges are non-critical.
\end{prop}

\proof
Recall from the proof of Theorem~\ref{thm:SG4struct} that the double cover mapping the drum $D_{k+1,2k+2}$ to the Schrijver graph $\SG(2k+2,k)$ maps the vertex $(i,j)$ of the drum to the vertex $V$ of the Scrijver graph satisfying that the two edges in $C_n\setminus V$ are the edge connecting $i+j-1$ to $i+j$ and the one connecting $j-i$ to $j-i+1$ (all calculations here and in the rest of this proof are modulo $2k+2$). Here $V$ is regular if and only if these edges are non-adjacent, so it is non-regular exactly if $i+j=j-i$ or $i+j-1=j-i+1$. The former equality happens exactly if $i=k+1$, while the latter happens exactly if $i=1$. Thus, the vertices in the top and bottom layer are mapped to non-regular vertices, while the rest of the layers are mapped to regular vertices. The isomorphism between the reduced drum and the Schrijver graph is obtained as the factor of the above map, therefore the vertices in the bottom layer of the reduced drum (corresponding to identification of a vertex from the bottom layer of the drum and its opposite vertex in the top layer) correspond to non-regular vertices in $\SG(2k+2,k)$, while the rest of the vertices in the reduced drum correspond to regular vertices as claimed.

Theorem~\ref{characterize} characterizes the critical edges in the Schrijver graph $\SG(2k+2,k)$ as the interlacing edges. By Proposition~\ref{regular}, all edges incident to a regular vertex are interlacing, so all edges of the reduced drum $D'_{k+1,2k+2}$ are critical except perhaps some edges in the bottom layer. Using the argument from the previous paragraph we can see that the vertex $(1,j)$ of the reduced drum corresponds to vertex $V_j=\{a\in[2k+2]\setminus\{j\}\mid a-j\hbox{ even}\}$ of the Schrijver graph. A neighbor of $(1,j)$ in the bottom layer of the reduced drum is $(1,i)$ with $i-j$ odd. If $j-1,j+1\in V_i$, then these nodes of $C_n$ are not separated by the set $V_j$, so the edge $V_iV_j$ is not interlacing making this edge (and the edge $(1,i)(1,j)$ of the reduced drum) non-critical. We have $j-1,j+1\in V_i$ unless $i=j-1$ or $i=j+1$. In these two cases the edge $V_iV_j$ is indeed interlacing, making it (and the edge $(1,i)(1,j)$ of the reduced drum) critical. \qed

\subsection{Critical edges of $\SG(n,2)$}\label{otnegy}

\begin{defi}
The \emph{length} $l(\{i,j\})$ of a vertex $\{i,j\}$ in the Schrijver graph $\SG(n,2)$ is the distance of $i$ from $j$ in the graph $C_n$, that is $l(\{i,j\})=\min(|i-j|,n-|i-j|)$. We have $2\le l(\{i,j\})\le\lfloor n/2\rfloor$.
\end{defi}

The following theorem characterizes the critical edges in $\SG(n,2)$. It shows that all interlacing edges are critical for $n\le7$ but this is not true for $n\ge8$, further most interlacing edges are non-critical as $n$ grows. We remark that this is consistent with a new result of Kaiser and Stehl\'{\i}k \cite{KSlect} which proves the edge-color-criticality of a spanning subgraph of $\SG(n,2)$ they obtain by deleting several non-interlacing edges.

\begin{thm}\label{notcharacterize}
An edge $XY$ of the Schrijver graph $\SG(n,2)$ is critical if and only if it is interlacing and either $l(X)\le3$ or $l(Y)\le3$.
\end{thm}

\proof
For brevity we write $\SG$ for $\SG(n,2)$.

For the if part consider an edge $XY$ of $\SG$ satisfying the conditions in the theorem. If both $X$ and $Y$ has length $2$, then $C_n[X\cup Y]$ is a path on four vertices, while if one of them has length $2$, while the other is longer, then $C_n[X\cup Y]$ consists of a path on three vertices and an isolated vertex. Proposition~\ref{basic} claims that the edge $XY$ is critical in $\SG$ in both cases.

In the remaining case either $X$ or $Y$ has length $3$. By symmetry we may assume that $X=\{1,4\}$ and $Y=\{2,i\}$ for some $i>4$. We define a coloring $f$ of the vertices of the graph $\SG$ as follows. If a vertex $U$ has an element $j>4$, $j\ne i$, then we set $f(U)$ to such a value $j$. This uses $n-5$ colors and creates no monochromatic edges. All vertices not colored yet are subsets of $\{1,2,3,4,i\}$. We use two special colors for these vertices, namely we set $f(U)=0$ if $U=\{1,3\}$, $U=\{1,i\}$ or $U=\{3,i\}$ and set $f(U)=-1$ if $U=X$, $U=Y$, $U=\{2,4\}$ or $U=\{4,i\}$. Note that this covers all remaining vertices and $\{1,i\}$ and $\{4,i\}$ may not even be vertices of the Schrijver graph depending on the value of $i$.

Clearly, the color class $0$ forms an independent set in the Schrijver graph while the color class $-1$ induces the single monochromatic edge $XY$. In total we use $n-3$ colors, so the coloring shows that the edge $XY$ is critical in $\SG$ as claimed.

For the only if part of the theorem let us consider a critical edge $XY$ of $\SG$. By Theorem~\ref{noncritical} $XY$ is interlacing, so we may assume $X=\{i,k\}$, $Y=\{j,l\}$ with $1\le i<j<k<l\le n$. Our goal is to prove that either $X$ or $Y$ has length at most $3$.

As $XY$ is critical we have a proper $(n-3)$-coloring $f$ of $\SG\setminus\{XY\}$. We identify vertices of $\SG$ with the corresponding edge of the complement $\overline C_n$ of $C_n$. This makes $f$ an edge-coloring of $\overline C_n$ where non-adjacent edges receive distinct colors except in the case of $X$ and $Y$. As $f$ cannot be a proper coloring of $\SG$, $X$ and $Y$ must indeed receive the same color. We call the color of $X$ and $Y$ red. Any further red edge must be adjacent to both $X$ and $Y$, so beyond $X$ and $Y$ only $ij$, $jk$, $kl$ and $li$ may be red (some of these may not even be edges in $\overline C_n$). Further, $ij$ and $kl$ are not adjacent, so at most one of them is red and similarly at most one of $jk$ and $li$ is red. Thus, we have at most $4$ red edges.

Any other color class contains pairwise adjacent edges, so these color classes
form either triangles of stars. Note that no color class can be a star
centered at $i$ or $k$, because then we could  modify $f$ by recoloring $X$
to this color to obtain a proper $(n-3)$-coloring of $SG$. Similarly no color
class can be a star centered at $j$ or $l$. Also, no two distinct color
classes can be stars with the same center as then we could join these classes
and recolor $X$ to a new color to obtain a proper $(n-3)$-coloring of $\SG$.

Let $s$ be the number of triangular color classes. We must have $n-4-s$ color
classes forming a star then. (Note that the class of the red edges is neither
triangular nor forms a star.)
We remove the center of these stars from $[n]$ to obtain a size $s+4$ subset
$Q$. In case a color class consists of a single edge we designate one end
vertex as the center and remove only that. Note that by the observation in the
previous paragraph $i,j,k,l\in Q$. Clearly, all edges of $\overline C_n[Q]$ belong to a triangular color class or are red. Out of the $s+4\choose2$ edges of the complete graph on $Q$ at most $s+4$ belong to the cycle $C_n$, the rest must be colored by the $s$ triangular color classes (at most $3$ edges per class) or must be red (at most $4$ edges), so we have
$${s+4\choose2}-(s+4)\le3s+4.$$

The inequality shows $s\le 2$. We use a case analysis to conclude that either $X$ or $Y$ has length at most $3$. We spell out this case analysis here for completeness, but the reader is probably better off proving this statement oneself than reading the rest of this proof.

We start with the case $s=0$. There is no triangular color class, so all edges in $\overline C_n[\{i,j,k,l\}]$ must be red. But as $ij$ and $kl$ cannot be both red, one of them must belong to $C_n$. Similarly, one of $jk$ and $li$ must also belong to $C_n$ and this makes the length of $X$ or $Y$ two.

In case $s=1$ we have $|Q|=5$. We have at most seven edges in two colors in
$\overline C_n[Q]$. This means that $C_n[Q]$ has at least $3$ edges. If it has
at least four edges, then at most one edge of $\overline C_n[Q]$ has length
larger than $3$, so one of $X$ and $Y$ has length at most $3$ and we are
done. So we assume that $C_n[Q]$ has exactly $3$ edges. This means that
$\overline C_n[Q]$ has $7$ edges, and among these a triangle is one color
class and $4$ red edges form another. Now consider the node $m\in
Q\setminus\{i,j,k,l\}$. The edges incident to $m$ are not red, therefore $2$
of the $4$ edges connecting $m$ to $i$, $j$, $k$ and $l$ must come from the
triangular color class, the other two from $C_n$. Say
$\{i,j,k,l\}=\{a,b,c,d\}$ with $ma$ and $mb$ from the triangular color class
and $mc$ and $md$ from $C_n$. Now $ab$ must also come from from the triangular
color class, therefore the third edge of $C_n[Q]$ must connect one of $a$ or
$b$ to one of $c$ or $d$. We may assume by symmetry that it is $bc$. As $bc$,
$cm$, and $md$ are all from $C_n$ the length of the edge $cd$ is $2$, while the length of $bd$ is $3$. One of these two edges is either $X$ or $Y$, so we are done again.

In case $s=2$ we have $|Q|=6$. Let $Q=\{i,j,k,l,m,m'\}$. Out of the edges
connecting $m$ to the other $5$ elements of $Q$ at most two belong to $C_n$,
none are red, so at least three belong to the two triangular color
classes. This implies that both the triangular color classes must include
edges incident to $m$. Similar claim can be said about $m'$. Thus, in particular, the pair $mm'$ should belong to both triangular color classes that is clearly impossible.
This finishes the proof of the theorem. \qed

\section{Relation to quadrangulations of nonorientable surfaces}\label{sect:quadr}

A graph is said to quadrangulate a surface if it can be embedded in the
surface in such a way that all the faces of the embedding are quadrangles.
The Schrijver graph ${\rm SG}(6,2)$ shown on Figure~\ref{fig:62} quadrangulates
the Klein bottle as shown on Figure~\ref{fig:quadr62}.

\begin{figure}[!htbp]
\centering
\includegraphics[scale=0.25]{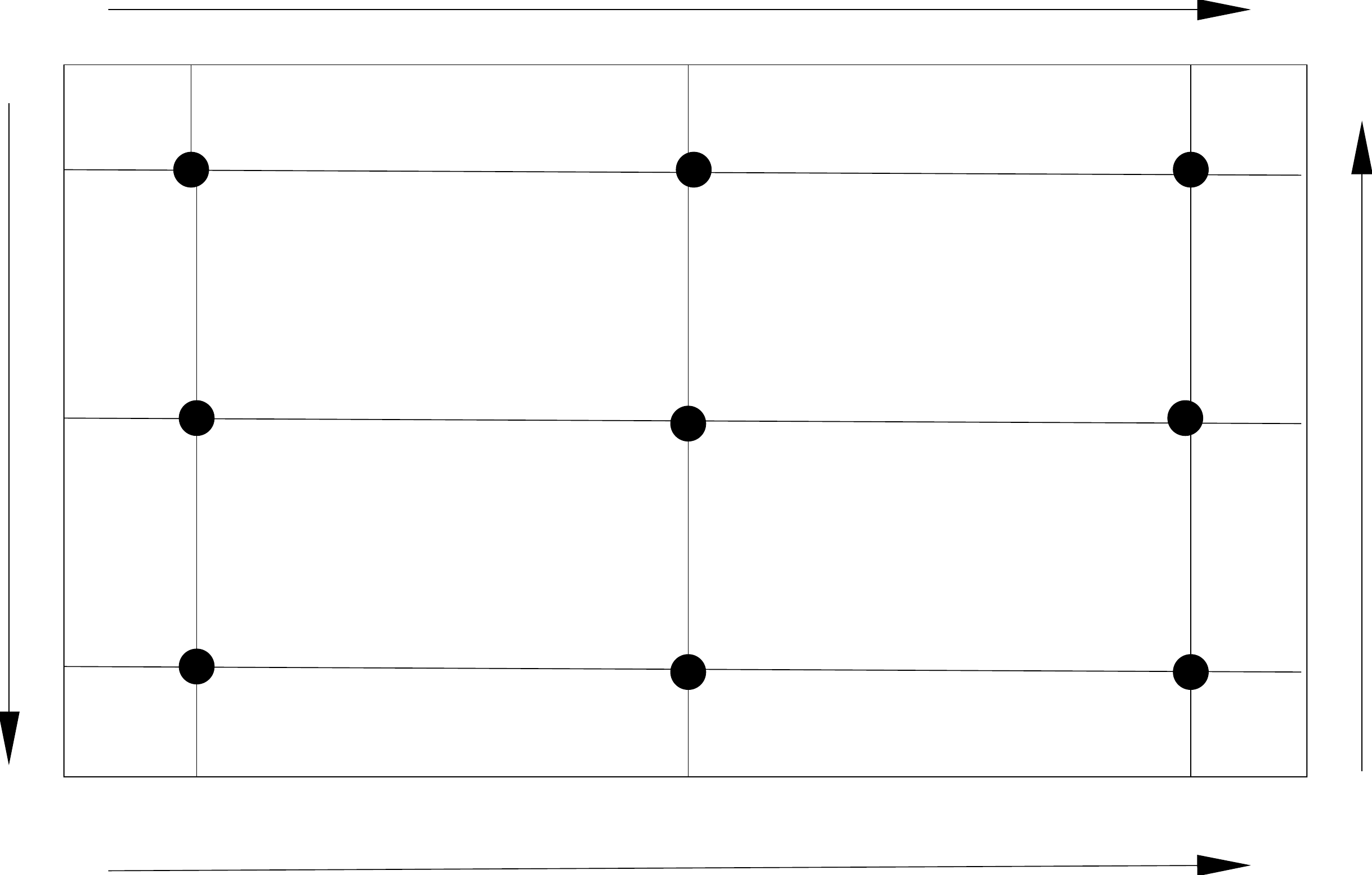}
\caption{The graph $Y_2={\rm SG}(6,2)$ drawn as a quadrangulation of the Klein
  bottle.}
\label{fig:quadr62}
\end{figure}

Similarly, we can draw $X_{h,n}$ (recall its definition from Section~\ref{sect:proof}
given before stating Theorem~\ref{bigyo}) as a quadrangulation of the Klein bottle by simply drawing a centrally symmetric grid of $2h-1$ rows and $n/2$ columns in the rectangle representing the Klein bottle. The vertices of the top and bottom rows induce the M\"obius ladder $M_n$, the vertices of the middle row induce $C_{n/2}$, while symmetric pairs of other rows induce a copy of $C_n$. Recall that for even $k$ the graph $X_{k/2+1,2k+2}$ is a spanning subgraph of $D'_{k+1,2k+2}=\SG(2k+2,k)$.
See Figure~\ref{fig:quadr146} for the drawing of $X_{4,14}$ which is a subgraph of $D'_{7,14}=SG(14,6)$ as a quadrangulation of the Klein bottle.

\begin{figure}[!htbp]
\centering
\includegraphics[scale=0.25]{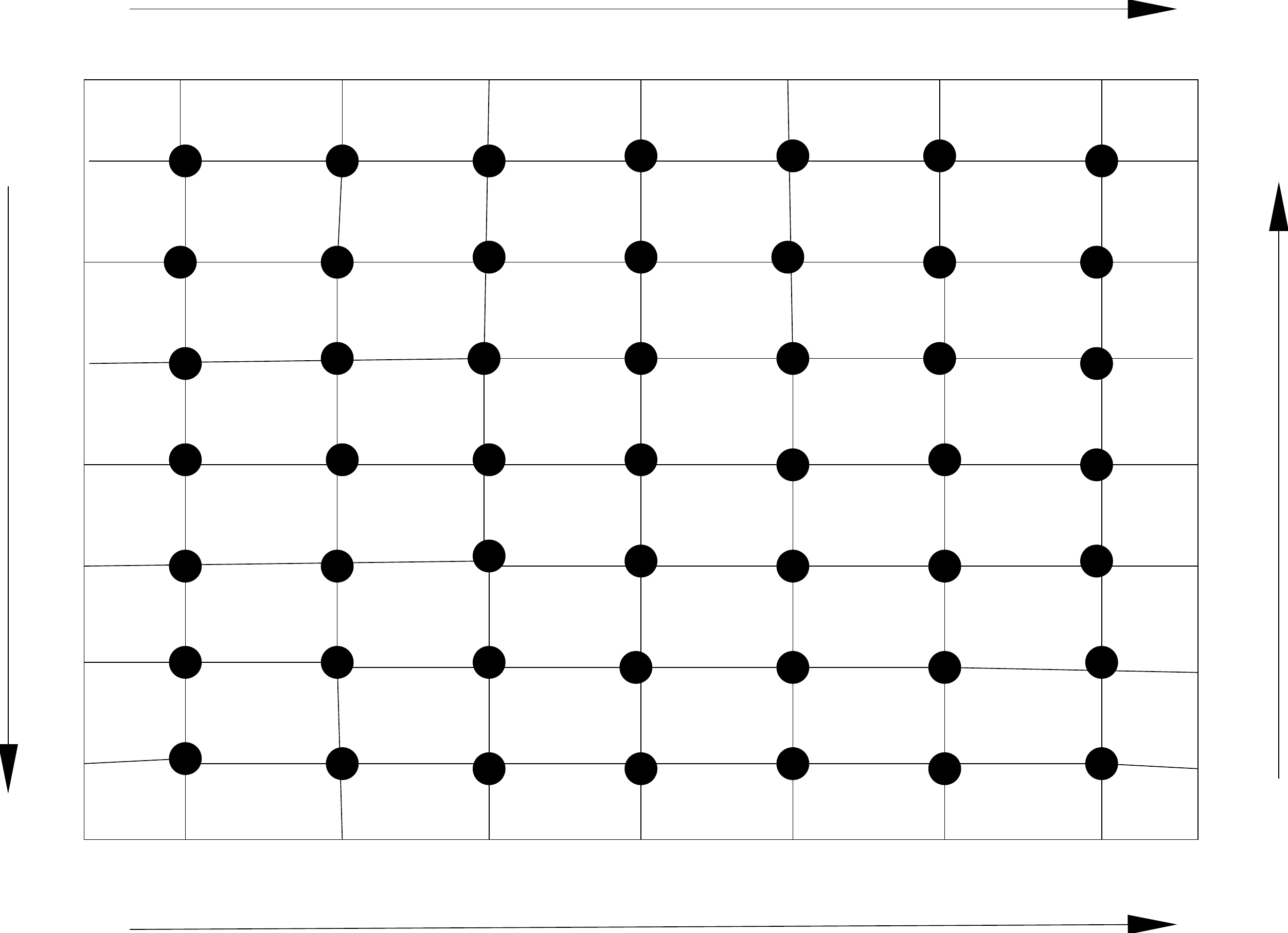}
\caption{The subgraph $Y_6=X_{4,14}$ of ${\rm SG}(14,6)$ drawn as a quadrangulation of the Klein
  bottle.}
\label{fig:quadr146}
\end{figure}

\medskip
\par\noindent
A spanning subgraph of ${\rm SG}(2k+2,k)$ quadrangulates the Klein bottle also in
case of odd $k>1$, but its drawing is a little more tricky. The
smallest case, ${\rm SG}(8,3)$ can be seen on Figure~\ref{fig:quadr83}.

\begin{figure}[!htbp]
\centering
\includegraphics[scale=0.25]{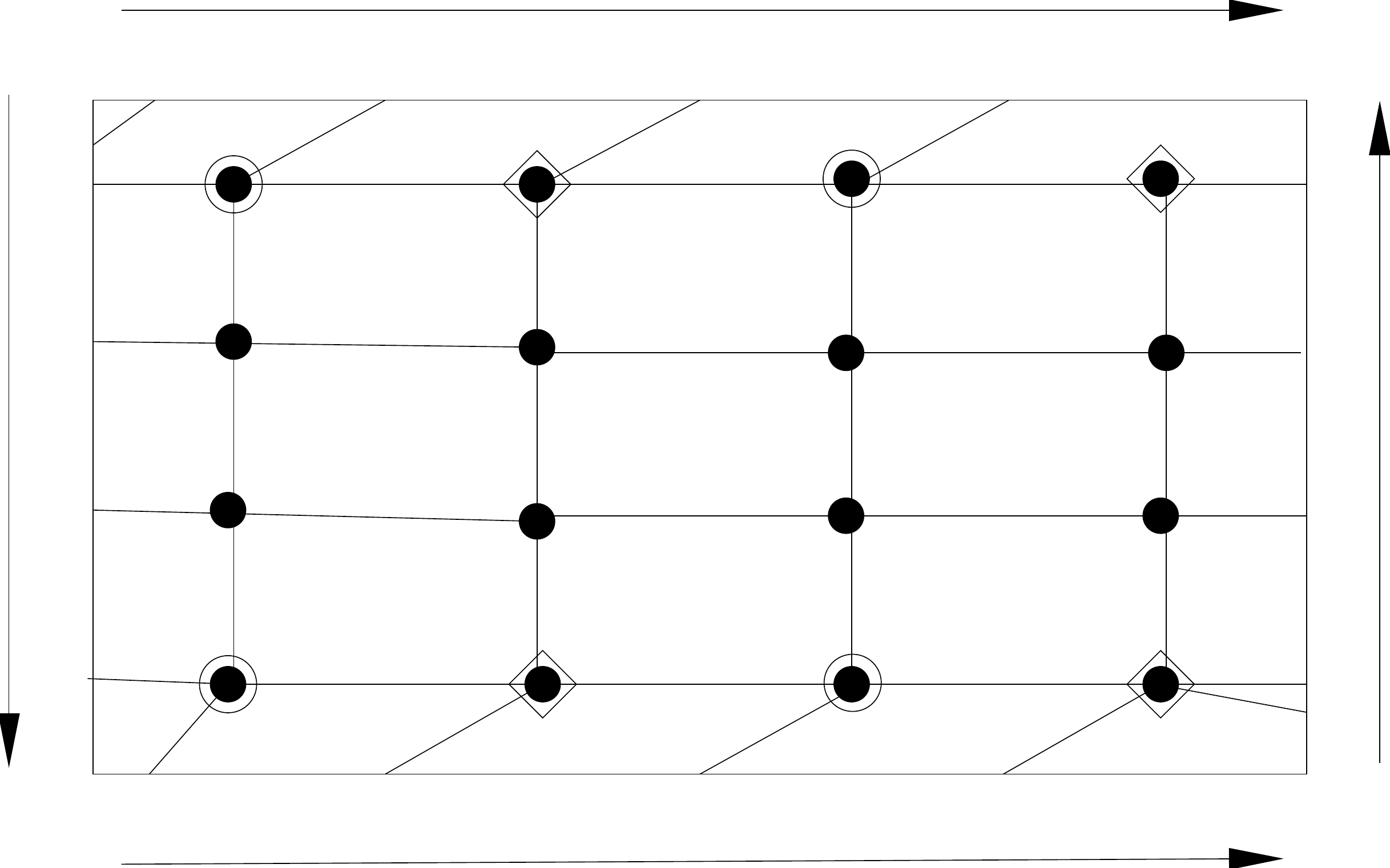}
\caption{The subgraph $Y_3$ of ${\rm SG}(8,3)$ drawn as a quadrangulation of the Klein
  bottle. The circled and squared vertices belong to the ``outer cycle''
  vertices of Figure~\ref{fig:83} that form a complete bipartite graph
  there. The circles and squares refer to their different color classes in
  that complete bipartite graph. Here
  only $12$ edges run among them out of the original $16$ ones, four edges of ${\rm SG}(8,3)$ are deleted to form this quadrangulation.}
\label{fig:quadr83}
\end{figure}

\medskip
\par\noindent
Again, it is not hard to see how to obtain a similar subgraph of
${\rm SG}(2k+2,k)$ quadrangulating the Klein bottle for any odd $k\ge 3$.
One has to draw the same grid-like picture with $k+1$ horizontal and vertical lines, where the top and bottom segment of the vertical
lines is made skew, so that the corresponding edge connects the top vertex of a vertical line with the bottom vertex of the next vertical line to the right. The top vertex of the rightmost line is not connected this way to any vertex (it has degree $3$), while there is an extra edge connecting the bottom left and bottom right vertices as in
Figure~\ref{fig:quadr83}. We will denote this graph by $Y_k$. In case $k$ is even, we write $Y_k$ to denote $X_{k/2+1,2k+2}$.

Note that $Y_k$ is typically not the only subgraphs of ${\rm SG}(2k+2,k)$ quadrangulating the Klein bottle. We chose it for its relative simplicity.

\smallskip
\par\noindent
Concerning the relevance of quadrangulations of the Klein bottle we remark that a closely related graph $G$ was demonstrated by Kaiser, Stehl\'{\i}k, and \v{S}krekovski \cite{KSS} to provide a counterexample to
a conjecture about square-free monomial ideals in polynomial rings. Graph $G$, originally defined by Gallai \cite{Gall} as remarked in \cite{KSS}, is simply the graph we obtain by putting one more vertical line into the picture of $\SG(6,2)$ on Figure~\ref{fig:quadr62}. For some other relations between algebraic questions and surface quadrangulations cf. \cite{MST}.

\medskip
\par\noindent
We can obtain a quadrangulation of the projective plane similarly. If we place
the $k+1$ by $k+1$ grid symmetrically in the rectangle representing the
projective plane the opposite corner vertices are connected by two parallel
edges, a horizontal one and a vertical one. Removing one edge from each of these pairs of parallel edges (say the vertical one) we obtain the quadrangulation $Z_k$ of the projective plane. Just as $Y_k$ is a subgraph of $\SG(2k+2,k)$, $Z_k$ is also a subgraph of $\SG(2k+2,k)$ (with the vertices located similarly in $Y_k$ and $Z_k$ in the two pictures representing the same vertex of $\SG(2k+2,k)$). Only the edges connecting the top and bottom row change. Note that the horizontal edges on these vertices form a cycle $C_{2k+2}$ in both cases. If $k$ is even the vertical edges connect opposite vertices of this cycle in $Y_k$. In case $k$ is odd the not-quite-vertical extra edges connect almost opposite vertices in the cycle. In both cases the subgraph formed by the cycle and these connecting edges is bipartite. In the case of $Z_k$ the vertical edges give diagonals extending the cycle $C_{2k+2}$ to a graph isomorphic to $P_{k+1}\Box P_2$, which is always bipartite.

Let us add a remark here on the appearance that  $Z_k$ is a subgraph of $Y_k$ at least in the case of even $k$. Indeed, one can delete two edges from the (picture of) $Y_k$ if $k$ is even to obtain the (picture of) $Z_k$. But this view is misleading. The two pictures may be equal but they are drawn in rectangles representing different surfaces and the change in the surface results in similarly placed segments joining different pairs of vertices. A similar relationship holds between $Y_k$ and $Z_k$ even if $k$ is odd if we draw the edges of $Z_k$ crossing the top edge of the bounding rectangle slanted as in Figure~\ref{fig:projquadr83}. This change of drawing does not effect the graph drawn (it is still $Z_k$ quadrangulating the projective plane), but now the picture is contained in the usual picture representing $Y_k$ (although in this case we cannot delete a few edges in $Y_k$ to obtain the exact same picture). To emphasize that the view of considering $Z_k$ as part of $Y_k$ is wrong we note that though the graph $Z_2$ happens to be a subgraph of $Y_2={\rm SG}(6,2)$, but the similar relation is not true for $Z_3$ and $Y_3$ (shown on Figures~\ref{fig:projquadr83} and \ref{fig:quadr83}) nor for $Z_4$ and $Y_4$ or beyond (for example for $Z_6$ and $Y_6$ shown on Figures~\ref{fig:projquadr146} and \ref{fig:quadr146}).

\begin{figure}[!htbp]
\centering
\includegraphics[scale=0.20]{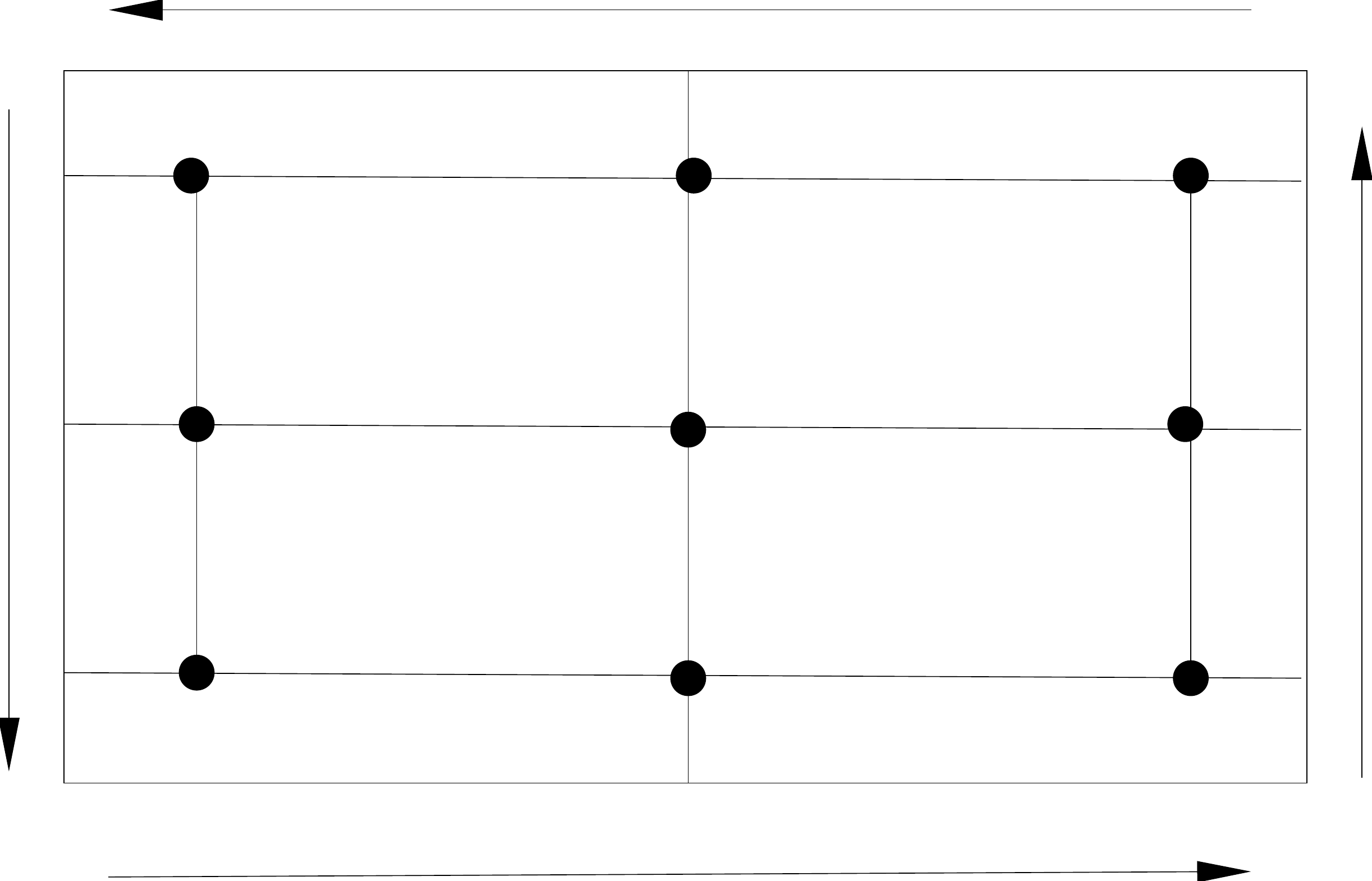}
\caption{The subgraph $Z_2$ of ${\rm SG}(6,2)$ drawn as a quadrangulation of the projective plane.}
\label{fig:projquadr62}
\end{figure}

\begin{figure}[!htbp]
\centering
\includegraphics[scale=0.20]{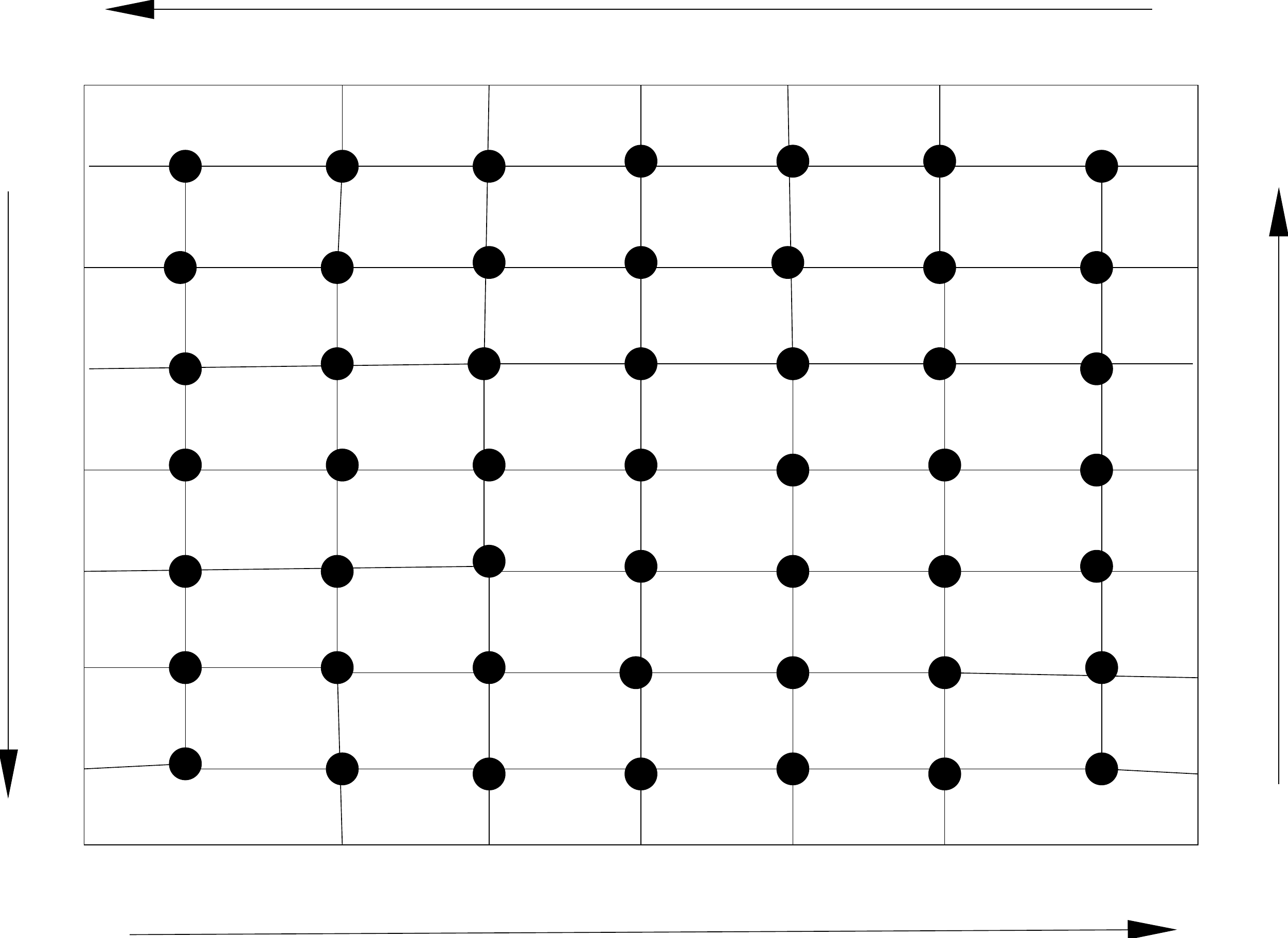}
\caption{The subgraph $Z_6$ of ${\rm SG}(14,6)$ drawn as a quadrangulation of the projective plane.}
\label{fig:projquadr146}
\end{figure}

\begin{figure}[!htbp]
\centering
\includegraphics[scale=0.25]{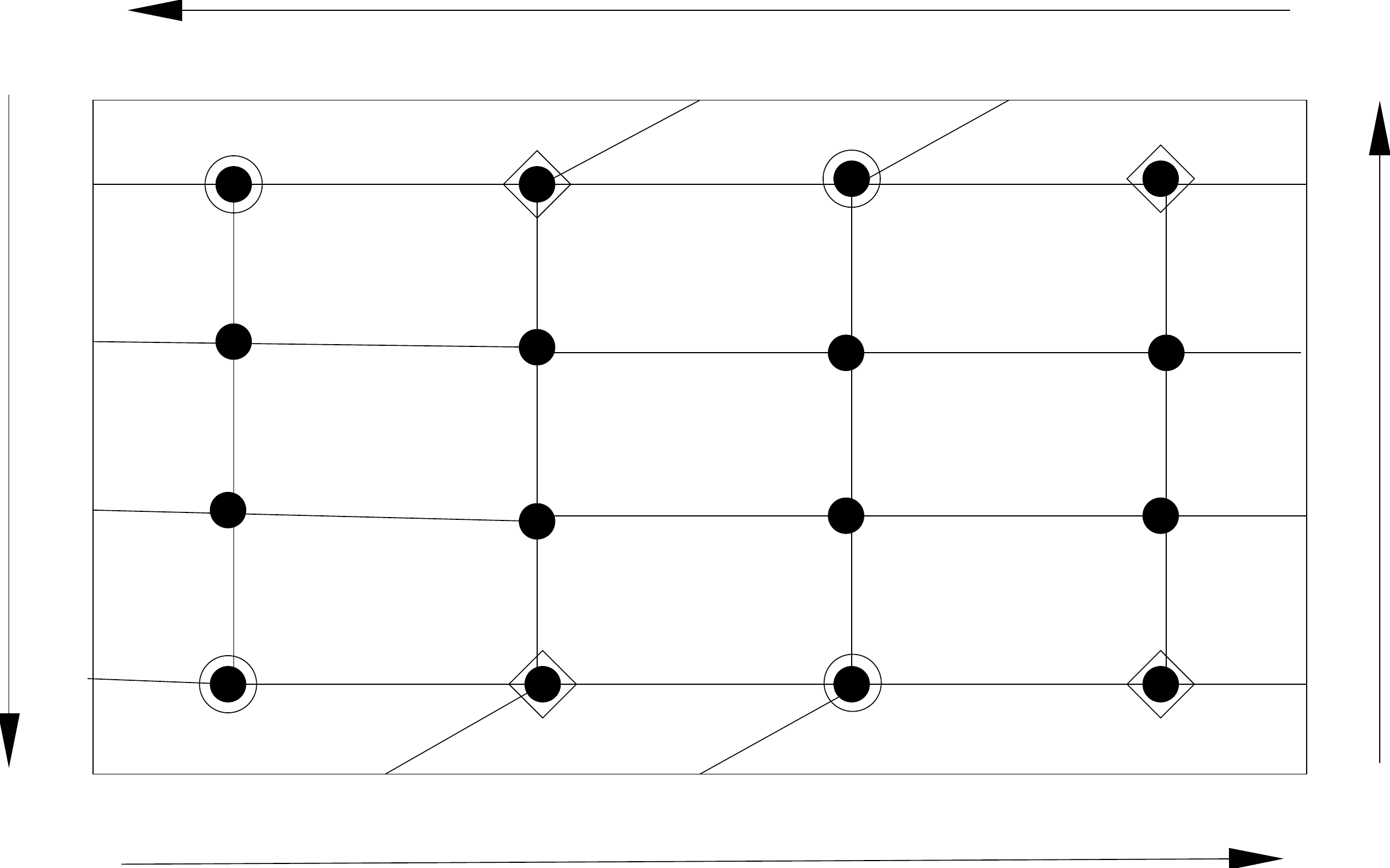}
\caption{The subgraph $Z_3$ of ${\rm SG}(8,3)$ drawn as a quadrangulation of the projective plane.}
\label{fig:projquadr83}
\end{figure}

\bigskip
\par\noindent
Chromatic properties of surface quadrangulations are widely investigated, see for example  \cite{Y, AHNNO, MS, MST}. Higher dimensional quadrangulations of
projective spaces are also investigated in \cite{KS1, KS2}. Youngs \cite{Y}
proved the surprising result that any quadrangulation of the projective plane
is either bipartite or $4$-chromatic, that is, it simply cannot have chromatic
number $3$. Archdeacon, Hutchinson, Nakamoto, Negami, Ota \cite{AHNNO} and
independently Mohar and Seymour \cite{MS} generalized this result to what are
called odd quadrangulations of nonorientable surfaces. They proved that if a
graph $G$ quadrangulates a nonorientable surface in such a way that cutting
the surface along an odd cycle of $G$ results in an orientable surface, then
the chromatic number of $G$ is at least $4$. Note that these results, along
with the observation above that a subgraph of every $4$-chromatic Schrijver
graph quadrangulates the Klein bottle, or the one that another subgraph quadrangulates the projective plane, already imply that these graphs, and thus also the graphs
${\rm SG}(2k+2,k)$ are not $3$-colorable. Youngs' theorem and the projective quadrangulations give this as $Z_k$ is not bipartite, in case of quadrangulating the Klein bottle one also has to verify the property that cutting along the edges of some odd cycle the surface would become orientable. It is not hard to find such a cycle though in our pictures. (In Figures~\ref{fig:quadr62}~and~\ref{fig:quadr146} any odd cycle formed by some ``vertical line'' would do, while only a little modification is needed in case of odd $k$.)
We also remark that the proof given in \cite{AHNNO} also makes use of the observation we stated in Lemma~\ref{lem:nonzero1}.

\medskip
\par\noindent
Kaiser and Stehl\'{\i}k \cite{KS1} define quadrangulations of topological spaces in any dimension as an extension of the quadrangulation of surfaces and generalize Youngs' theorem for quadrangulations of the projective space of higher dimensions.
In \cite{KS2} they prove that a spanning subgraph of every Schrijver graph
${\rm SG}(n,k)$ quadrangulates the projective space of dimension $n-2k$. (With
their generalization of Youngs' theorem this gives the right value for the
chromatic number of ${\rm SG}(n,k)$.) Their construction is quite involved and
even in the case $n=2k+2$ which is relevant for us, it is not easy to compare
it to our construction. (In fact, Figures~\ref{fig:projquadr83} and
\ref{fig:paltq83} demonstrate that the spanning subgraph of $\SG(2k+2,k)$
quadrangulating the projective plane need not be unique.)
They also conjecture that their spanning subgraph (for every $n$ and $k$) is edge-color-critical and note that in case of $n=2k+2$ it follows from a result of Gimbel and Thomassen \cite{GT}. Indeed, Gimbel and Thomassen proved that if $G$ is a graph in the projective plane such that all contractible cycles have length at least four, then $G$ is $3$-colorable if and only if $G$ does not contain a nonbipartite qadrangulation of the projective plane. This implies that a nonbipartite quadrangulation of the projective plane is edge-color-critical if all of its contractible four-cycles are faces.

\subsection{Critical edges in the Klein bottle quadrangulating subgraphs}\label{critklein}

As we have already mentioned at the end of the last section,
a non-bipartite quadrangulation of the projective plane is edge-color-critical if all contractible 4-cycles are facial by a
result of Gimbel and Thomassen~\cite{GT}. This applies to the quadrangulations $Z_k$. Here we consider the
subgraph $Y_k$ of $\SG(2k+2,k)$ that quadrangulates the Klein bottle.
We saw that $Z_2$ happens to be a proper subgraph of $Y_2$ but no similar relation holds between $Z_k$ and $Y_k$ for $k\ge3$. This implies that $Y_2$ is not edge-color-critical, but does not say anything about $Y_k$ for $k\ge3$.
We will see that the graph $Y_3$ (see on
Figure~\ref{fig:quadr83}) is \emph{not} edge-color-critical either but $Y_k$ is edge-color-critical for $k\ge4$.

To see the structure of $Y_k$ clearly we introduce the following graphs.

\begin{defi}\label{defi:Lk}
For even values of $k$ let $L_k$ be obtained from the cycle $C_{2k+2}$ by adding the \emph{extra edges} connecting opposite vertices in this cycle. So $L_k$ is the M\"obius ladder $M_{2k+2}$ in this case.
\smallskip\par\noindent
For odd values of $k$ let us obtain $L_k$ from the cycle $C_{2k+2}$ by adding the \emph{extra edges} $i(i+k)$ for $i\in[k+1]$.
\end{defi}

Note that the graph $L_k$ is bipartite for all values of $k$.

Considering $\SG(2k+2,k)$ as the reduced drum $D'_{k+1,2k+2}$ as described by Lemma~\ref{ize}, we can see (cf.\ the figures in the previous section) that $Y_k$ is isomorphic to the graph we obtain from the reduced drum $D'_{k+1,2k+2}$ if we change its bottom layer (a complete bipartite graph) to its subgraph isomorphic to $L_k$ in the natural way (i.e., such that the vertex $i$ of $L_k$ corresponds to the vertex $(1,i)$ of $D'_{k+1,2k+2}$). This means deleting some edges of $D'_{k+1,2k+2}$ from its bottom layer. We identify $L_k$ with this subgraph of $Y_k$ so we can refer to extra edges in $Y_k$: these are the extra edges of $L_k$ in the bottom layer.

Note that $Y_2$ is isomorphic to $\SG(6,2)$, therefore its critical edges are described by Theorem~\ref{characterize} (or, rather by Proposition~\ref{geomch}). It has $3$ non-critical edges. The following theorem deals with all larger values of $k$.

\begin{thm}\label{thm:Klcrit}
For $k\ge 4$ the subgraph $Y_k$ of $\SG(2k+2,k)$ quadrangulating the Klein bottle is edge-color-critical. In $Y_3$ all edges but two are critical.
\end{thm}

\proof
Set again $\SG=\SG(2k+2,k)$ for the sake of brevity of notation.
Having the characterization of critical edges of $\SG$ by Proposition~\ref{geomch}, we only have to prove that the non-critical edges of $\SG$ that still appear in $Y_k$ for $k\ge 4$ become critical there (while exactly two of them become critical for $k=3$). The edges of $Y_k$ non-critical in $\SG$ are exactly the extra edges of $L_k$.

Let us denote by $U_k$ the following graph.
$$V(U_k):=V(P_2\Box C_{2k+2})$$
and for even $k\ge 4$
$$E(U_k):=E(P_2\Box C_{2k+2})\cup \{(1,i)(1,j): \{i,j\}\in E(L_k)\},$$
while for odd $k\ge 3$
$$E(U_k):=E(P_2\Box C_{2k+2})\cup \{(1,i)(1,j): ij\in E(L_k)\}\cup\{(2,i)(2,i+k+1): i\in [k+1]\}.$$
Note that the bottom layer of $U_k$ induces $L_k$, while the top layer induces the cycle $C_{2k+2}$ or the M\"obius ladder $M_{2k+2}$ depending on the parity of $k$. When speaking of an extra edge of $U_k$ we mean an extra edge of $L_k$ in the bottom layer of $U_k$.

We observe that it will be enough to prove that deleting any extra edge of $U_k$ for odd $k\ge 5$, the remaining graph has a proper $3$-coloring, while deleting any extra edge of $U_k$ for even $k\ge 4$, the remaining graph has a proper $3$-coloring in such a way that the opposite vertices of the induced $C_{2k+2}$ of its second layer get identical colors. (The case of $k=3$ will be dealt with separately, but in fact we have $Y_3\cong U_3$.) This is because for odd $k\ge 5$ $Y_k$ differs from $U_k$ only by having some ``middle layers'' inducing cycles $C_{2k+2}$ between the two layers of $U_k$. However, once a $3$-coloring of $U_k\setminus\{e\}$ is found for some extra edge $e$ of $U_k$, then it can be extended to $Y_k$ in a very simple way. Indeed, copy such a $3$-coloring to the two bottommost layers of $Y_k$ and once the $j$th layer is colored, color the $(j+1)$th layer by the rule $c(j+1,i):=c(j,i+1)$. This will give a proper $3$-coloring of $Y_k$ when $k$ is odd. The even case is similar, the only difference is that the topmost layer of $Y_k$ is in fact a cycle $C_{k+1}$ of half length. But it is clear that its coloring is equivalent to coloring a double length cycle $C_{2k+2}$ in such a way that the opposite vertices get identical colors. (Identifying these opposite vertices will give the required coloring of the half-length cycle.)

We consider first the case when $k\ge 4$ is even. We have to prove that if any extra edge is deleted from $U_k$, then the remaining graph has a proper $3$-coloring $c$ for which $c(2,i)=c(2,i+k+1)$ for every $i\in [k+1]$. We may assume by symmetry that we remove the edge $(1,1)(1,k+2)$, since an automorphism of $U_k$ can bring this edge to any other extra edge of $U_k$ while permuting the second layer within itself. We give the coloring $c$ explicitly. Reading the colors of the vertices in the first layer in order we obtain
$$012(01)^{k/2-1}020(12)^{k/2-1},$$
while in the second layer we similarly obtain
$$201(20)^{k/2-1}201(20)^{k/2-1}.$$
Here the exponential notation $(ab)^i$ means that the colors $a$ and $b$ alternately repeat $i$ times.

Note that this coloring gives the same color (namely $0$) to both endpoints of the deleted edge $(1,1)(1,k+2)$, but creates no other monochromatic edge in $U_k$ and gives identical colors to opposite vertices of the cycle in the second layer. See Figure~\ref{fig:U6minus} for an illustration of this $3$-coloring.

\begin{figure}[!htbp]
\centering
\includegraphics[scale=0.35]{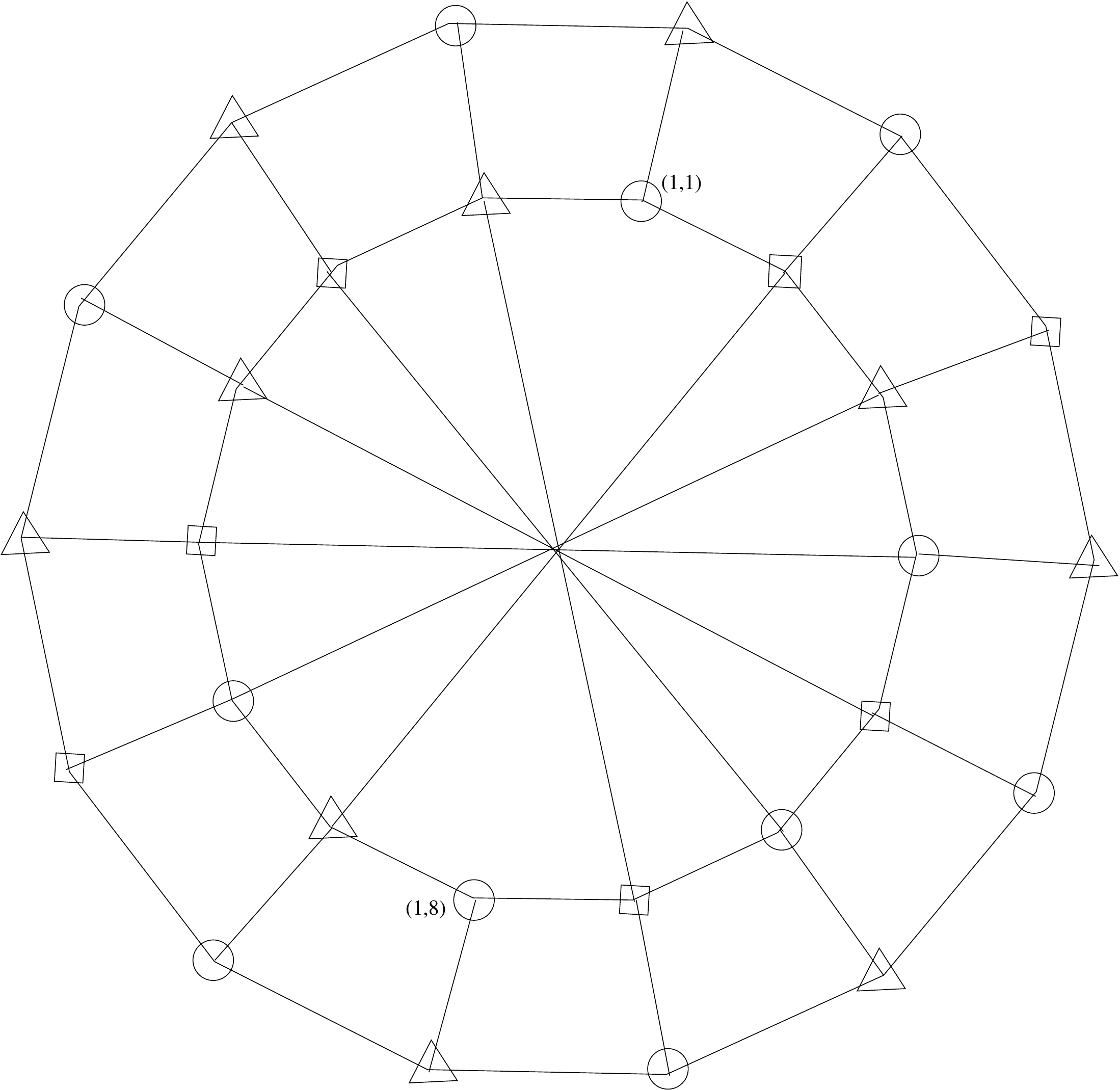}
\caption{Our $3$-coloring of $U_6\setminus\{e\}, e=(1,1)(1,8)$ where opposite
  vertices of the second layer get identical colors. Colors $0,1$, and $2$ are
indicated by circle, square, and triangle, respectively.}
\label{fig:U6minus}
\end{figure}

\medskip\par\noindent
Now assume $k$ is odd. We should give a proper $3$-coloring of the vertices of $U_k$ after any one extra edge is deleted. Now $U_k$ has only one non-trivial automorphism, so it is not enough to consider the removal of a single extra edge.


Let us first consider the removal of the extra edge $e=(1,i)(1,i+k)$ from $U_k$ for some $2\le i\le k$. Note that $e$ is also an extra edge in $U_{k-1}$, so we have a proper $3$-coloring $c$ of $U_{k-1}\setminus\{e\}$ that gives the same color to opposite vertices in the second layer. We extend $c$ to a proper 3-coloring of $U_k\setminus\{e\}$ as follows. We set $c(1,2k+1):=c(1,1)$ and $c(2,2k+1):=c(2,1)$ leaving the vertices $(1,2k+2)$ and $(2,2k+2)$ uncolored for now. We claim we have created no monochromatic edge yet. Indeed, the edges of $U_k$ among the vertices $(1,j)$ for $j\in[2k]$ are also edges of $U_{k-1}$, so (other than the deleted edge $e$) none is monochromatic. All neighbors of the vertex $(1,2k+1)$ in the first layer of $U_k$ are also neighbors of $(1,1)$ in $U_{k-1}$, so these edges did not become monochromatic either. Clearly, no ``vertical'' edge $(1,j)(2,j)$ became monochromatic either for $j\in[2k+1]$, nor did any of the edges $(2,j)(2,j+1)$ for $j\in[2k]$. Finally, we have to check the ``diagonal'' edges $(2,j)(2,j+k+1)$ for $j\in[k]$. these are not monochromatic, because $c(2,j)=c(2,j+k)\ne c(2,j+k+1)$ by the properties of our coloring of $U_{k-1}\setminus\{e\}$.

It remains to further extend $c$ to the two remaining vertices of $Y_k\setminus\{e\}$. We start with the vertex $(2,2k+2)$. It has three neighbors already colored, two of which are $(2,1)$ and $(2,2k+1)$ with identical colors, so we can find a color for the vertex $(2,2k+2)$ not represented in its neighborhood and thus keeping our coloring proper. Finally we consider the last remaining vertex $(1,2k+2)$. It has three neighbors in $U_k$ out of which $(1,1)$ and $(1,2k+1)$ have the same color, so we can extend our proper 3-coloring to this vertex too. This proves that $\chi(U_k\setminus\{e\})\le3$ as needed.
\medskip

We still have to consider the case when the deleted edge $e$ is one of $(1,1)(1,k+1)$ and $(1,k+1)(1,2k+1)$. These two cases are symmetric (an automorphism of $U_k$ takes one to the other), so we may assume that the deleted edge is the former one. We show an explicit proper $3$-coloring $c$ of $U_k\setminus\{e\}$ in this case. With the same notation as before, reading the colors of the vertices in the first layer in order we obtain
$$0(21)^{(k-1)/2}0(10)^{(k-1)/2}21,$$
while in the second layer we obtain
$$1(02)^{(k-1)/2}1(01)^{(k-1)/2}02.$$
This coloring gives color 0 to both endpoints of $e$, but makes no other edge of $U_k$ monochromatic. This makes $c$ a proper $3$-coloring of $U_k\setminus\{e\}$ and finishes the proof that $Y_k$ is edge-color critical for $k\ge4$.

Finally, we consider the case $k=3$, namely the graph $Y_3=U_3$. The $3$-coloring given explicitly in the previous paragraph works in the $k=3$ case, too, that is, it gives a proper $3$-coloring of $Y_3\setminus\{(1,1)(1,4)\}$ and thus it shows that the edge $(1,1)(1,4)$ is critical in $Y_3$. By symmetry, the edge $(1,4)(1,7)$ is also critical in $Y_3$. However, neither of the two remaining extra edges, namely $(1,2)(1,5)$ and $(1,3)(1,6)$ are critical. Even if we delete both of these edges, the remaining graph $Q$ is $4$-chromatic, albeit not isomorphic to $Z_3$. One way to see that $Q$ has no proper $3$-coloring is to notice that the remaining two extra edges, that is, $(1,1)(1,4)$ and $(1,4)(1,7)$ divide the cycle $C_8$ in the bottom layer into three quadrangles and then applying the argument of Lemmata~\ref{lem:nonzero1}, \ref{lem:ketkor}, and \ref{lem:nonzero2} on the two
layers of $Q$.
Alternatively, one can show that $Q$ quadrangulats the projective plane (as $Z_3$ does also) and (since it is not bipartite) we have $\chi(Q)=4$ by Youngs' theorem. This completes the proof of the theorem.
\hfill$\Box$

\begin{figure}[!htbp]
\centering
\includegraphics[scale=0.25]{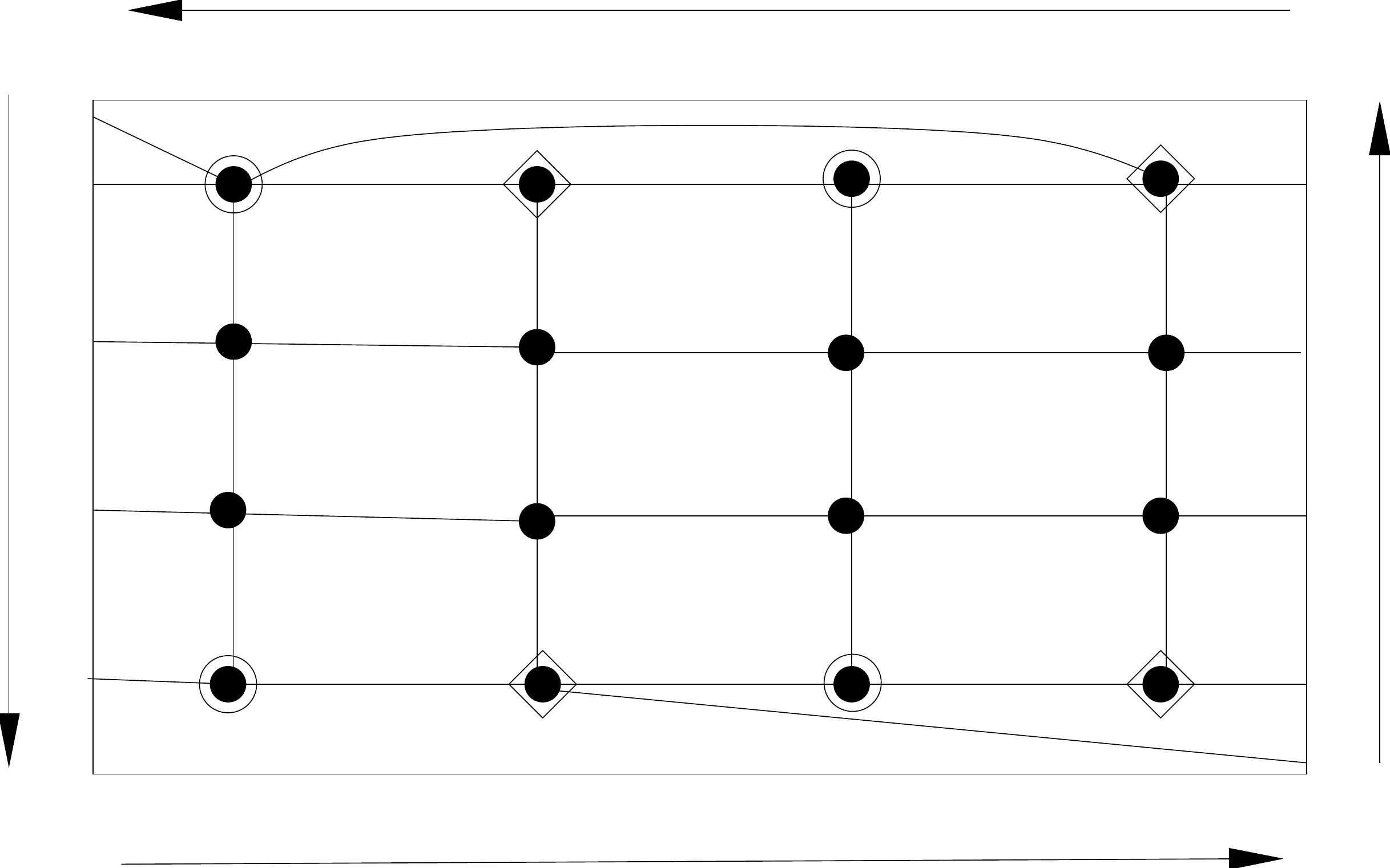}
\caption{A subgraph of $\SG(8,3)$ quadrangulating the projective plane
  different from the one on Figure~\ref{fig:projquadr83}. (Note that unlike
  the graph on Figure~\ref{fig:projquadr83} this one has a vertex of degree
  $5$.) This is the quadrangulation mentioned at
  the end of the proof of Theorem~\ref{thm:Klcrit}.}
\label{fig:paltq83}
\end{figure}

\section{Open problems}\label{open}

In this final section we come back to the question of critical edges of Schrijver graphs of higher chromatic numbers. Recall that our results characterized them for $4$-chromatic Schrijver graphs as being the interlacing edges (Theorem~\ref{characterize}), but this characterization does not hold in general (see Theorem~\ref{notcharacterize}). We cannot even come up with a conjecture for the full characterization, nevertheless, we formulate an open problem (that, if true, would generalize Theorem~\ref{2regular}) and a conjecture that would characterize the critical edges for large Schrijver graphs with small chromatic number.

\begin{problem}
Are all edges incident to a regular vertex in a Schrijver graph critical?
\end{problem}

\begin{conj}
For any $d$ and large enough $k$ (that is $k>k_d$ for some threshold $k_d$ depending on $d$) all the interlacing edges of the Schrijver graph $\SG(2k+d,k)$ are critical.
\end{conj}

\section{Acknowledgements}

We thank Tom\'a\v{s} Kaiser and Mat\v{e}j Stehl\'{\i}k for a useful
conversation and in particular for drawing our attention to the paper
\cite{LPSV}.

\end{document}